\documentclass[12pt]{amsart}
\usepackage{amsmath}
\usepackage{amssymb}
\usepackage{amsmath,amssymb,amsthm,amscd}
\usepackage{relsize}
\usepackage{caption}

\usepackage{txfonts}
\usepackage{amsbsy}
\usepackage{graphicx,color,epstopdf}

\numberwithin{equation}{section}

\newtheorem{prop}{Proposition}[section]
\newtheorem{theo}{Theorem}[section]
\newtheorem{lemm}{Lemma}[section]

\def\begeq{\begin{equation}}
\def\endeq{\end{equation}}

\begin{document}

\title{On the Planar $L_p$-Minkowski Problem}
\author{Shi-Zhong Du}
\thanks{The author is partially supported by STU Scientific Research Foundation for Talents (SRFT-NTF16006), and Natural Science Foundation of Guangdong Province (2019A1515010605)}
  \address{The Department of Mathematics,
            Shantou University, Shantou, 515063, P. R. China.} \email{szdu@stu.edu.cn}

\renewcommand{\subjclassname}{%
  \textup{2010} Mathematics Subject Classification}
\date{Sep. 2020}

\begin{abstract}
   In this paper, we study the planar $L_p$-Minkowski problem
      \begin{equation}\label{e0.1}
        u_{\theta\theta}+u=fu^{p-1}, \ \ \theta\in{\mathbb{S}}^1
      \end{equation}
   for all $p\in{\mathbb{R}}$, which was introduced by Lutwak \cite{L2}. A detailed exploration of \eqref{e0.1} on solvability will be presented. More precisely, we will prove that for $p\in(0,2)$, there exists a positive function $f\in C^\alpha({\mathbb{S}}^1), \alpha\in(0,1)$ such that \eqref{e0.1} admits a nonnegative solution vanishes somewhere on ${\mathbb{S}}^1$. In case $p\in(-1,0]$, a surprising \textit{a-priori} upper/lower bound for solution was established, which implies the existence of positive classical solution to each positive function $f\in C^\alpha({\mathbb{S}}^1)$. When $p\in(-2,-1]$, the existence of some special positive classical solution has already been known using the Blaschke-Santalo inequality \cite{CW}. Upon the final case $p\leq-2$, we show that there exist some positive functions $f\in C^\alpha({\mathbb{S}}^1)$ such that \eqref{e0.1} admits no solution. Our results clarify and improve largely the planar version of Chou-Wang's existence theorem \cite{CW} for $p<2$. At the end of this paper, some new uniqueness results will also be shown.
\end{abstract}

\maketitle\markboth{$L_p$ Minkowski Problem}{$L_p$-surface measure}

\tableofcontents

\section{Introduction}

The Minkowski problem is to determine a convex body with prescribed curvature or other similar geometric data. It plays a central role in the theory of convex bodies. Various Minkowski problems \cite{A,CL,CY,GLL,GM,HLYZ,P} have been studied especially after Lutwak \cite{L1,L2}, who proposes two variants of the Brunn-Minkowski theory including the dual Brunn-Minkowski theory and the  $L_p$ Brunn-Minkowski theory. Besides, there are singular cases such as the logarithmic Minkowski problem and the centro-affine Minkowski problem \cite{BLYZ,CW}.

The main purpose of this paper is to study the $L_p$-Minkowski problem for different exponents $p$ in the plane. Given a convex body $\Omega$ in ${\mathbb{R}}^{n+1}$ containing origin, for each $x\in{\mathbb{S}}^n$, let $r(x)$ be a point on $\partial\Omega$ whose outer unit normal is $x$. The support function $h(x)$ is defined to be
   $$
    h(x)\equiv{r}(x)\cdot x, \ \ \forall x\in{\mathbb{S}}^n.
   $$
For a uniformly convex $C^2$ body, the matrix $[\nabla^2_{ij}h+h\delta_{ij}]$ is positive definite, where $[\nabla^2_{ij}h]$ stands for the Hessian tensor of $h$ acting on an orthonormal frame of ${\mathbb{S}}^n$. Conversely, any $C^2$ function $h$ satisfying $[\nabla^2_{ij}h+h\delta_{ij}]>0$ determines a uniformly convex $C^2$ body $\Omega_h$. Direct calculation shows the standard surface measure of $\partial\Omega$ is given by
    $$
     dS\equiv\frac{1}{n+1}\det(\nabla^2_{ij}h+h\delta_{ij})d {\mathcal{H}}^{n}|_{{\mathbb{S}}^n}
    $$
for $n$-dimensional Hausdorff measure $d {\mathcal{H}}^{n}$.
 It's well known that classical Minkowski problem looks for a convex body such that its standard surface measure matches a given Radon measure on ${\mathbb{S}}^n$. In \cite{L2}, Lutwak introduces the $L_p-$surface measure $dS_p\equiv h^{1-p}dS$ on $\partial\Omega$. The corresponding $L_p$-Minkowski problem is to look for a convex body whose $L_p-$surface measure is equal to a prescribed function.

 Parallel to the classical Minkowski problem, the $L_p$-Minkowski problem boils down to solve the fully nonlinear equation
   \begin{equation}\label{e1.1}
       \det(\nabla^2_{ij}h+h\delta_{ij})=fh^{p-1}, \ \ \forall x\in{\mathbb{S}}^n
   \end{equation}
in the smooth category. In their corner stone paper \cite{CW}, Chou and Wang study the $L_p$-Minkowski problem and obtain the following result.

 \begin{theo}\label{t1.1}
   Considering the $L_p$-Minkowski problem \eqref{e1.1}, one has

 (1) When $p>n+1$, there exists a unique positive solution in $C^{2,\alpha}({\mathbb{S}}^n)$ for each positive function $f\in C^\alpha({\mathbb{S}}^n)$. And

 (2) when $p=n+1$, there exists a unique pair $(h,\lambda)$ for $0<h\in C^{2+\alpha}({\mathbb{S}}^n)$ and $0<\lambda\in{\mathbb{R}}$ satisfying
   \begin{equation}\label{e1.2}
     \det(\nabla^2_{ij}h+h\delta_{ij})=\lambda fh^n
   \end{equation}
 for each positive function $f\in C^\alpha({\mathbb{S}}^n), \alpha\in(0,1)$. And

 (3) when $1<p<n+1$, \eqref{e1.1} has a generalized non-negative solution in the sense of Aleksandrov for each $f\in L^\infty({\mathbb{S}}^n), f\geq f_0$, where $f_0$ is some positive constant. And

 (4) when $p\in(-n-1,1)$, \eqref{e1.1} has a generalized non-negative solution in the sense of Aleksandrov for each $f\in L^\infty({\mathbb{S}}^n), f\geq f_0$, where $f_0$ is some positive constant. Moreover, if $p\in(-n-1,-n+1]$ and $f\in C^\alpha({\mathbb{S}}^n)$ for some $\alpha\in(0,1)$, this special solution is positive and in $C^{2,\alpha}({\mathbb{S}}^n)$.
 \end{theo}

In the existence theorem of Chou-Wang \cite{CW}, only the case $p>n+1$ was solved completely for all dimensions $n$ due to the validity of the maximum principle in this case. Since the lacking of a \textit{a-priori} positive lower bound to the solution of \eqref{e1.1}, the remaining cases $-n-1\leq p<n+1$ were analysed only partially on weak sense and the case $p<-n-1$ has not been discussed yet.

 For the planar case, the $L_p$-Minkowski equation \eqref{e1.1} becomes a second-order nonlinear ordinary differential equation
   \begin{equation}\label{e1.3}
     h_{\theta\theta}+h=fh^{p-1}, \ \ \mbox{ on } {\mathbb{S}}^1,
   \end{equation}
 where $\theta$ is the arc-length parameter of ${\mathbb{S}}^1$. As usually, we assume that $f$ is a positive solution belonging to $C^\alpha({\mathbb{S}}^1)$ for some $\alpha\in(0,1)$ and name the positive solution $h\in C^{2,\alpha}({\mathbb{S}}^1)$ to be the classical one.
 The existence of positive classical solution for $p>2$ has been obtained in \cite{CW}, meanwhile an eigenvalue version of \eqref{e1.3} was solved for $p=2$ there. When $p\in(0,2)$, we will prove \textit{a-priori} upper and lower bounds for the width function of the convex set, and show that there can not be an \textit{a-priori} positive lower bound for the positive classical solution in general.

  \begin{theo}\label{t1.2}
    For $p\in(0,2)$ and positive function $f\in C^\alpha({\mathbb{S}}^1)$, there exists a positive constant $C_{p,f}\geq1$ depending only on $p$ and $f$, such that
      \begin{equation}\label{e1.5}
         C_{p,f}^{-1}\leq w_\Omega^-\leq w_\Omega^+\leq C_{p,f}, \ \ \mbox{ on } {\mathbb{S}}^1,
      \end{equation}
where
     \begin{eqnarray*}
       w_\Omega^-&\equiv& \min_{\theta\in[0,2\pi)}(u(\theta)+u(\theta+\pi)),\\
       w_\Omega^+&\equiv& \max_{\theta\in[0,2\pi)}(u(\theta)+u(\theta+\pi))
     \end{eqnarray*}
are minimal and maximal width of $\Omega$ respectively. However, there is some positive function $f\in C^\alpha({\mathbb{S}}^1)$ such that \eqref{e1.3} admits a nonnegative but not positive solution.
  \end{theo}

A proof of theorem \ref{t1.2} will be presented in Section 2 and Section 3, which was inspired by work of Chen-Li \cite{CL} for dual Minkowski problem. Although \textit{a-priori} positive lower bound for positive classical solution of \eqref{e1.3} can not be expected for $p\in(0,2)$, we still have the following surprising result for $p\in(-1,0]$.

 \begin{theo}\label{t1.3}
   Assuming $p\in(-1,0]$ and $f\in C^\alpha({\mathbb{S}}^1)$ is positive, there exists a positive constant $C_{p,f}\geq1$ depending only on $p, \min f$ and $||f||_{C^\alpha({\mathbb{S}}^n)}$, such that
      \begin{equation}\label{e1.8}
        C_{p,f}^{-1}\leq h\leq C_{p,f}, \ \ \mbox{ on } {\mathbb{S}}^1
      \end{equation}
    holds for all positive classical solution $h$ of \eqref{e1.3}. Consequently, for each positive function $f\in C^\alpha({\mathbb{S}}^1)$, there exists at least one positive classical solution to \eqref{e1.3}.
 \end{theo}

Comparing to the result of Chou-Wang \cite{CW} upon the planar setting, we have derived a new \textit{a-priori} positive lower bound \eqref{e1.8} to solution of \eqref{e1.3} in case $p\in(-1,0]$. As a result, one can obtain the existence of positive classical solution in stead of nonnegative weak solution in \cite{CW}. The proof to Theorem \ref{t1.3} will be presented in Section 4-6. For the case $p\in(-2,-1]$, the solvability of nonnegative solutions was already shown in \cite{CW} using Blaschke-Santalo's inequality and the method of variational. At the remaining range $p\leq-2$, we will prove the following non-existence result.

\begin{theo}\label{t1.4}
  Assuming $p=-2$, there exist some positive functions $f\in C^\alpha({\mathbb{S}}^1), \alpha\in(0,1)$ such that (1.3) is not solvable. When $p<-2$, a same result holds for some H\"{o}lder functions $f$ which is positive outside two polar of ${\mathbb{S}}^1$.
\end{theo}

The equation \eqref{e1.1} is invariant under all projective transformations on the sphere ${\mathbb{S}}^n$ in the case $p=-n-1$. Using the view point of centroaffine geometry for this Minkowski problem, Chou and Wang \cite{CW} found a striking necessary condition for solvability of \eqref{e1.1} for $p=-n-1$ in terms of derivative of $f$ along the projective vector field $\xi\in{\mathbb{S}}^n$. However, owing to the dependent of unknown solution $h$ and absence of positivity of kernel of $f$, this condition can not be applied directly to produce an explicit example of insolvability to our planar cases $p\leq-2$. Fortunately, using a new trigonometric relation for solution of \eqref{e1.3}, we have constructed explicit examples of $f$ ensuring insolvability of \eqref{e1.3}, even in the deeply negative case $p<-2$. We will prove Theorem \ref{t1.4} in Section 7.

Next we discuss the uniqueness of positive solutions to \eqref{e1.3}. Using the maximum principle, Chou-Wang have shown uniqueness result for $p>n+1$ in \cite{CW} for all dimension. The uniqueness result certainly can not be expected for $p=n+1$ due to the homogeneity of equation. When $p<n+1$, the uniqueness problem is much more subtle since lack of the maximum principle. There are only some partial results are known on the past. Chow has shown in \cite{C} for all $n\geq1$ and $p=1-n$, the uniqueness holds true for constant function $f\equiv1$. Using the invertible result for linearized equation, Lutwak showed in \cite{L2} for $n\geq1, p>1$, uniqueness holds for some special symmetric $f$. Later, Dohmen-Giga \cite{DG} and Gage \cite{G} have extended the result to $n=1, p=0$ for some symmetric function $f$. Contrary to the results in \cite{DG,G}, Yagisita \cite{Y} showed a surprising non-uniqueness result for $n=1, p=0$ and non-symmetric function $f$. Subsequently, Andrews showed uniqueness for $n=2, p=0$ and arbitrary positive function $f$. When $p\in(-n-1,-n-1+\sigma), 0<\sigma\ll1$, a counter example of uniqueness have also obtained by Chou-Wang in \cite{CW}. More recently, Jian-Lu-Wang \cite{JLW} have proven that for $p\in(-n-1,0)$, there exists at least a smooth positive function $f$ such that \eqref{e1.1} admits two different solutions. While a partial uniqueness result was established by Chen-Huang-Li-Liu \cite{CHLL} on origin symmetric convex bodies for $p\in(p_0,1), p_0\in(0,1)$, using the $L_p$-Brunn-Minkowski inequality. For the deep negative case $p\leq-n-1$, the situations are more complicated. As it is well known that for $n\geq1, p=-n-1$, all ellipsoids with the volume of the unit ball are all solutions of \eqref{e1.1} with $f\equiv1$, the uniqueness fails in this case. When $p<-n-1$, Andrews showed in \cite{An} that the uniqueness property is much more delicate even in dimension $n=1$.

 For its importance and difficulty, the issue of uniqueness of solutions has attracted much attention. The problems have been conjectured for a number of special cases, including in particular the case $f\equiv1, p=0$ by Firey \cite{F}, and the case $f\equiv1, p\in(-n-1,1)$ by Lutwak-Yang-Zhang \cite{HLX,LYZ}. This is the main purpose for us to discuss the problem. We will consider $p<2$ to planar equation \eqref{e1.3} and prove the following result.

\begin{theo}\label{t1.5}
Letting $f\equiv1$, we have

(1) If $p\in(-2,1)\cup(1,2)$, constant solution $h\equiv1$ is the unique positive classical solution of \eqref{e1.3}.

(2) If $p\in(-\infty,-7)$, there exist at least
 $$
   c_p\equiv\big[\sqrt{2-p}\big]_*-1\geq2
 $$
positive classical solutions to \eqref{e1.3}.

 (3) When $p=0$, the constant solution $h\equiv1$ is the unique positive classical solution to \eqref{e1.3}. While for $p=1, -2$ and each given $h_{min}\in(0,1)$, there exists a positive classical solution $h$ of \eqref{e1.3} satisfying
   \begin{equation}\label{e1.9}
      \min_{{\mathbb{S}}^1}h=h(0)=h_{min}.
   \end{equation}

 (4) For any $p<2, p\not=0, 1, -2$, there exists a positive constant $\sigma_p$ such that there is no positive classical solution $h$ of \eqref{e1.3} satisfying
    \begin{equation}\label{e1.10}
       1-\sigma_p<h_{min}<1.
    \end{equation}
\end{theo}

\vspace{10pt}

 The uniqueness result in part (1) was shown by Lutwak in \cite{L2} for $p\in(1,2)$, and shown by Chou-Zhu \cite{CZ} for $p\in(-2,1)$. We present here a different proof in case of $p\in[1/2,1)\cup(1,2)$. While, the nonuniqueness result in part (2) for $p<-7$ was new so far as we known. We will give the proof of Theorem \ref{t1.5} in Section 8. Moreover, a purely algebraic sufficient condition on $p$ ensuring the uniqueness property was given in Proposition \ref{p8.3}. Combining our uniqueness results with previously known ones, one has the following theorem for general positive function $f$.

  \begin{theo}\label{t1.6}
    Considering \eqref{e1.3} for each positive function $f\in C^\alpha({\mathbb{S}}^1)$, the follows hold:

  (1) For $p>2$, uniqueness holds for all $f$.\\

  (2) For $p\in(1,2)$, uniqueness holds for some special symmetric $f$.\\

  (3) For $p=0$, uniqueness holds for some special symmetric $f$ and fails to hold for some other non-symmetric $f$.\\

  (4) For $p\in(-2,0)$, uniqueness fails to hold for some $f$.\\

  (5) For $p\in(-\infty,-7)$, uniqueness fails to hold for $f\equiv1$.\\

  (6) Uniqueness fails to hold for $p=1,-2$ and $f\equiv1$.\\
  \end{theo}

Complete uniqueness result for $p>2$ was obtained by Chou-Wang in \cite{CW}, while the partial uniqueness result for $p>1$ was shown by Lutwak in \cite{L2}. When $p=0$, the uniqueness for symmetric $f$ has been proven in \cite{DG,G}, and the counterexample for non-symmetric $f$ was due to \cite{Y}. If $p\in(-2,0)$, the nonuniqueness for some special $f$ was given by \cite{JLW}. The non-unique result for $p\in(-\infty,-7)$ and constant function $f\equiv1$ was given in part (2) of our Theorem \ref{t1.5}. The final part (6) was known already owing to the explicitly examples as mention above.

\vspace{40pt}

\section{``Good shape" estimation for $0<p<2$}

This section is devoted to the case $0<p<2$, in which an \textit{a-priori} upper/lower bound for width function will be given. Roughly speaking, we will show that the convex body is of ``good shape" in the sense of Theorem \ref{t2.1}. Given a convex body $\Omega\subset{\mathbb{R}}^2$ containing origin, we denote its support function by
   $$
    u(\theta)\equiv\sup\big\{y\cdot x|\ y\in\Omega\big\}, \ \ \forall x=\left(\begin{array}{c}
            \cos\theta\\
            \sin\theta
         \end{array}
         \right)\in{\mathbb{S}}^1.
   $$
Conversely, letting $u$ be a positive function satisfying
    $$
     u_{\theta\theta}+u>0, \ \ \forall\theta\in{\mathbb{S}}^1,
    $$
we denote $\Omega_u$ to be the convex body determined by support function $u$. One has the expansion formula
   \begin{equation}\label{e2.1}
     {r}(\theta)=u\left(
         \begin{array}{c}
            \cos\theta\\
            \sin\theta
         \end{array}
         \right)+u_\theta\left(
         \begin{array}{c}
            -\sin\theta\\
            \cos\theta
         \end{array}
         \right)
   \end{equation}
for the radial boundary vector ${r}(\theta)\in\partial\Omega$, whose unit outer normal is given by $x=\left(
         \begin{array}{c}
            \cos\theta\\
            \sin\theta
         \end{array}
         \right)$. Consequently, the radial length function of $\Omega$ is defined by
  $$
   {l(\theta)}\equiv|{r}(\theta)|=\sqrt{u^2+u_\theta^2}.
  $$
Let's start with several elementary lemmas which would be used later.

\begin{lemm}\label{l2.1}
   Supposing that the solution $u$ of \eqref{e1.3}  is monotone increase on interval $(\theta_1,\theta_2)$, one has
      \begin{equation}\label{e2.2}
       \begin{cases}
        l^2(\theta_2)-l^2(\theta_1)\leq\frac{2f_+}{p}\Big(u^p(\theta_2)-u^p(\theta_1)\Big)\\
        l^2(\theta_2)-l^2(\theta_1)\geq\frac{2f_-}{p}\Big(u^p(\theta_2)-u^p(\theta_1)\Big)
       \end{cases}
      \end{equation}
  in case $p\not=0$ and
      \begin{equation}\label{e2.3}
       \begin{cases}
       \displaystyle l^2(\theta_2)-l^2(\theta_1)\leq2f_+\ln\frac{u(\theta_2)}{u(\theta_1)}\\[10pt]
      \displaystyle  l^2(\theta_2)-l^2(\theta_1)\geq2f_-\ln\frac{u(\theta_2)}{u(\theta_1)}
       \end{cases}
      \end{equation}
 in case $p=0$,  where $f_+, f_-$ are maximum and minimum of $f$ respectively. As a result, there hold
      \begin{equation}\label{e2.4}
       \begin{cases}
         u_{max}^2-\frac{2f_+}{p}u_{max}^p\leq u_{min}^2-\frac{2f_+}{p}u_{min}^p,\\
         u_{max}^2-\frac{2f_-}{p}u_{max}^p\geq u_{min}^2-\frac{2f_-}{p}u_{min}^p
       \end{cases}
      \end{equation}
 in case $p<2, p\not=0$ and
      \begin{equation}\label{e2.5}
       \begin{cases}
         u_{max}^2-2f_+\ln u_{max}\leq u_{min}^2-2f_+\ln u_{min},\\
         u_{max}^2-2f_-\ln u_{max}^p\geq u_{min}^2-2f_-\ln u_{min}
       \end{cases}
      \end{equation}
 in case $p=0$, where $u_{max}$ and $u_{min}$ are maximum and minimum of $u$ respectively. When $u$ is monotone decreasing function on $\theta\in(\theta_1,\theta_2)$, the above inequalities \eqref{e2.2} and \eqref{e2.3} will be reversed.
\end{lemm}

\noindent\textbf{Proof.} The proof to this lemma is elementary. In fact, multiplying \eqref{e1.3} by $u_\theta$, one has
   $$
     \Big(u^2+u_\theta^2\Big)_\theta\leq\frac{2f_+}{p}(u^p)_\theta
   $$
on the increasing arc $[\theta_1,\theta_2]$ and so obtains first inequality of \eqref{e2.2}. The second inequality of \eqref{e2.2} and \eqref{e2.3} are similarly. To show the first inequality of \eqref{e2.4}, one needs only to draw a picture for the function $F(u)\equiv u^2-\frac{2f_+}{p}u^p$ in case of $p<2, p\not=0$. Since there is only one minimal point of $F$ on ${\mathbb{R}}^+$ and $\lim_{u\to+\infty}F(u)=+\infty$, one can prove that no mater $u$ is monotone increasing from $u_{min}$ to $u_{max}$ or not, there always holds first inequality of \eqref{e2.4}. In fact, let's suppose that $u$ increases from minimal point $\vartheta_0$ to a first critical point $\vartheta_1>\vartheta_0$. If $\vartheta_1$ is exactly the maximal point of $u$, then first inequality of \eqref{e2.4} follows from first inequality of \eqref{e2.2} by replacing
   $$
    \theta_1=\vartheta_0, \ \ \theta_2=\vartheta_1.
   $$
If $\vartheta_1$ is only a local maximal point of $u$, one also has
   \begin{equation}\label{e2.6}
     u^2(\vartheta_1)-\frac{2f_+}{p}u^p(\vartheta_1)\leq u_{min}^2-\frac{2f_+}{p}u^p_{min}.
   \end{equation}
Now, assuming that $u$ decreases from $\vartheta_1$ to the first local minimal point $\vartheta_2>\vartheta_1$, after using $u(\vartheta_2)\geq u_{min}$ and picture of $F$, we still have
   \begin{equation}\label{e2.7}
     u^2(\vartheta_2)-\frac{2f_+}{p}u^p(\vartheta_2)\leq u_{min}^2-\frac{2f_+}{p}u^p_{min}.
   \end{equation}
By a bootstrap argument, the first inequality of \eqref{e2.4} was shown. The validity of second inequality of \eqref{e2.4} and \eqref{e2.5} can be verified similarly. The proof of the lemma was done. $\Box$\\

A second lemma follows from the maximum principle.

\begin{lemm}\label{l2.2}
  Supposing that $u$ is a positive classical solution of \eqref{e1.3} for $p<2$, one has
    \begin{equation}\label{e2.8}
      u_{max}\geq f_-^{\frac{1}{2-p}}, \ \ u_{min}\leq f_+^{\frac{1}{2-p}}.
    \end{equation}
 Furthermore, if $p\in(0,2)$, there holds
    \begin{equation}\label{e2.9}
      u_{max}\leq C_{p,f_+}
    \end{equation}
 for some positive constant $C_{p,f_+}$ depending only on $p$ and $f_+$.
\end{lemm}

\noindent\textbf{Proof.} \eqref{e2.8} is a direct consequence of the maximum principle of \eqref{e1.3}, and \eqref{e2.9} follows from \eqref{e2.4} together with Young's inequality. $\Box$\\

Now, given a convex body $\Omega$ with support function $u$, we define its width function by
   $$
     w_\Omega(\theta)\equiv u(\theta)+u(-\theta), \ \ \mbox{ on } {\mathbb{S}}^1
   $$
and denote
   $$
     w_\Omega^+\equiv\max_{{\mathbb{S}}^1}w_\Omega, \ \ w_\Omega^-\equiv\min_{{\mathbb{S}}^1}w_\Omega
   $$
to be its maximal and minimal width. By John's lemma \cite{J}, there exists an ellipsoid $E$ such that
   $$
      E-\xi\subset\Omega-\xi\subset 2(E-\xi).
   $$
After a rotation if necessary, one may assume that $a\geq b>0$ and
  $$
    E\equiv\Bigg\{(z_1,z_2)\in{\mathbb{R}}^2\Big|\ \frac{(z_1-\xi_1)^2}{a^2}+\frac{(z_2-\xi_2)^2}{b^2}=1\Bigg\},
  $$
where $\xi=(\xi_1,\xi_2)\in{\mathbb{R}}^2$ is the center of $E$. Thus, there clearly holds that
  \begin{equation}\label{e2.10}
   \begin{cases}
     w_{\Omega}^+/2\leq u_{max}\leq w_{\Omega}^+,\ \  u_{min}\leq w_{\Omega}^-\\
     a=w_E^+\leq w_\Omega^+\leq 2w_E^+=2a, \ \ b=w_E^-\leq w_\Omega^-\leq 2w_E^-=2b
   \end{cases}
  \end{equation}
by comparison. Furthermore, by an inversion if necessary, one could assume
  $$
    u(\pi)\leq u(0),\ \ u(3\pi/2)\leq u(\pi/2)
  $$
and set
   $$
    L\equiv u(0), \ \ l\equiv u(\pi/2), \ \ d\equiv u(\pi).
   $$
There will be an equivalent property
   \begin{equation}\label{e2.11}
       a/2\leq L\leq 2a, \ \ \frac{4b}{35}\leq l\leq 2b,
   \end{equation}
where the inequality $l\geq \frac{4b}{35}$ follows from a comparison of the intersections of lines $z_1=\frac{4}{5}L$ and $z_1=\xi_1$ with $\Omega$ (see Figure1). In fact, suppose that the line $z_1=L$ intersects with $\partial\Omega$ at $A$ and the line $z_1=\frac{4}{5}L$ intersect with $\partial\Omega$ at $B, C$. Drawing two rays $\overrightarrow{AB}$ and $\overrightarrow{AC}$ intersect with the line $z_1=\xi_1$ at $B'$ and $C'$. Since
   $$
     -\frac{3L}{4}\leq\xi_1\leq \frac{3L}{4},
   $$
we have
   $$
     \frac{35l}{2}\geq\frac{7L/4}{L/5}|BC|\geq|B'C'|\geq2b.
   $$
So, one gets $l\geq\frac{4b}{35}$.

\begin{figure}
  \includegraphics[width=5.1in,height=2.55in]{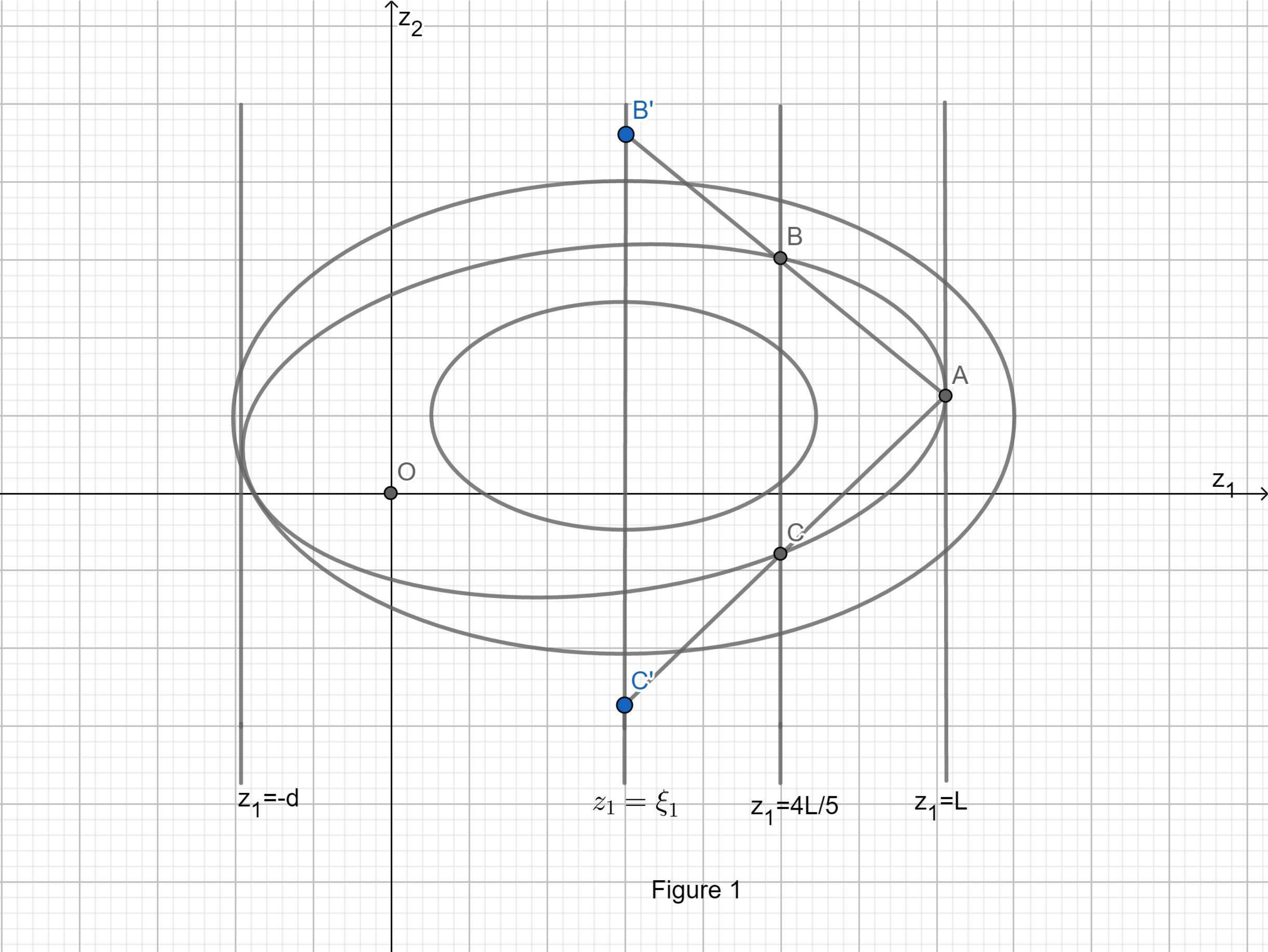}
  \caption{Estimation on lower bound of $l$}
\end{figure}

Our first main result is the following estimation to upper and lower bounds of width length function. In other words, we will prove the convex body $\Omega$ would be a ``round" one rather than a ``thin" one, and would be neither ``small" nor ``large" in the scale.

\begin{theo}\label{t2.1}
   For each $p\in(0,2)$ and positive function $f\in C^\alpha({\mathbb{S}}^1)$, there exists a positive constant $C_{p,f}$ depending only on $p$ and $f_\pm$, such that
     \begin{equation}\label{e2.12}
       C_{p,f}^{-1}\leq w_{\Omega_u}^-\leq w_{\Omega_u}^+\leq C_{p,f}
     \end{equation}
   holds for positive classical solution $u$ of \eqref{e1.3}.
\end{theo}

Set $m^*$ to be the $z_1$ coordinate of ${r}\Big(\frac{\pi}{2}\Big)$ and define $\theta_0-\theta_1$ by
  $$
    {r}(\theta_0)=(0,\varphi(0)), \ \ {r}(\theta_1)=(4L/5,\varphi(4L/5)),
  $$
where $\varphi(s), s\in[-d,L]$ is the graph function of upper component of boundary $\partial\Omega$. Before proving the theorem, let's first quote a geometric lemma by Chen-Li \cite{CL}.

\begin{lemm}\label{l2.3}
  When $m^*>0$, the function ${l(\theta)}$ and $u(\theta)$ are both monotone decreasing functions from $\pi/2$ to $\theta_0$. If $m^*\leq\frac{3L}{4}$, the function ${l(\theta)}$ and $u(\theta)$ are also both monotone decreasing functions from $\theta_1$ to $\pi/2$, as long as
    \begin{equation}\label{e2.13}
       w_{\Omega_u}^+/w_{\Omega_u}^-\geq100.
    \end{equation}
\end{lemm}

\noindent\textbf{Proof.} The first part in case $m^*\geq\frac{L}{2}$ follows from monotone increasing of $\varphi(s)$ on interval $s\in[0,m^*]$. To show the second part in case $m^*\leq\frac{L}{2}$, one needs only to use the bound
   $$
    |\varphi'(s)|\leq1/10, \ \ \forall s\in[m^*,4L/5]
   $$
thanks to the assumption \eqref{e2.13}. $\Box$\\

\noindent\textbf{Proof of Theorem \ref{t2.1}.} Suppose on the contrary, one may assume that \eqref{e2.13} holds. When $m^*\geq \frac{L}{2}$, by Lemma \ref{l2.3} and \ref{l2.1}, we have
   \begin{equation}\label{e2.14}
     l^2(\pi/2)-l^2(\theta_0)\leq\frac{2f_+}{p}\Big(u^p(\pi/2)-u^p(\theta_0)\Big)\leq\frac{2f_+}{p}l^p.
   \end{equation}
Noting that upon \eqref{e2.13}, there holds also
   \begin{equation}\label{e2.15}
     l^2(\pi/2)-l^2(\theta_0)\geq L^2/4-l^2\geq L^2/5.
   \end{equation}
Another hand, it yields from \eqref{e2.10} and \eqref{e2.11} that
   \begin{equation}\label{e2.16}
      \frac{w_{\Omega_u}^+}{4}\leq L\leq 2w_{\Omega_u}^+.
   \end{equation}
Combining \eqref{e2.14}-\eqref{e2.16} with
    \begin{equation}\label{e2.17}
      \frac{w_{\Omega_u}^+}{2}\leq u_{max}\leq w_{\Omega_u}^+
    \end{equation}
and Lemma \ref{l2.2}, one concludes that
    \begin{equation}\label{e2.18}
      C_{p,f}^{-1}\leq l\leq L\leq C_{p,f}.
    \end{equation}
Thus, a contradiction holds with our assumption by using \eqref{e2.10}, \eqref{e2.11} and \eqref{e2.18}.

If $m^*\leq3L/4$ and \eqref{e2.13} holds, we use
  \begin{equation}\label{e2.19}
    l^2(\theta_1)-l^2(\pi/2)\leq\frac{2f_+}{p}\Big(u^p(\theta_1)-u^p(\pi/2)\Big)\leq Cl^p
  \end{equation}
to replace \eqref{e2.14} and use
   \begin{equation}\label{e2.20}
     l^2(\theta_1)-l^2(\pi/2)\geq C^{-1}L^2
   \end{equation}
to replace \eqref{e2.15}, where
   $$
    \begin{cases}
      \tan\alpha\equiv\Big|\varphi'\Big(\frac{4L}{5}\Big)\Big|\leq\frac{l}{L/5},\\
      \tan\beta\equiv\frac{\varphi(4L/5)}{4L/5}\leq\frac{l}{4L/5}
    \end{cases}
   $$
and
   \begin{eqnarray*}
     u(\theta_1)&=&\sqrt{\Big(\frac{4L}{5}\Big)^2+\varphi^2\Big(\frac{4L}{5}\Big)}\sin(\alpha+\beta)\\
       &\leq&\sqrt{(4L/5)^2+l^2}\sin(\alpha+\beta)\leq CL\cdot\frac{l}{L}=Cl
   \end{eqnarray*}
have been used in second inequality of \eqref{e2.19} (see Figure 2). A similar argument as above gives the positive lower bound of $l$ and thus \eqref{e2.12} was drawn. $\Box$\\

\begin{figure}
  \includegraphics[width=5.1in,height=2.55in]{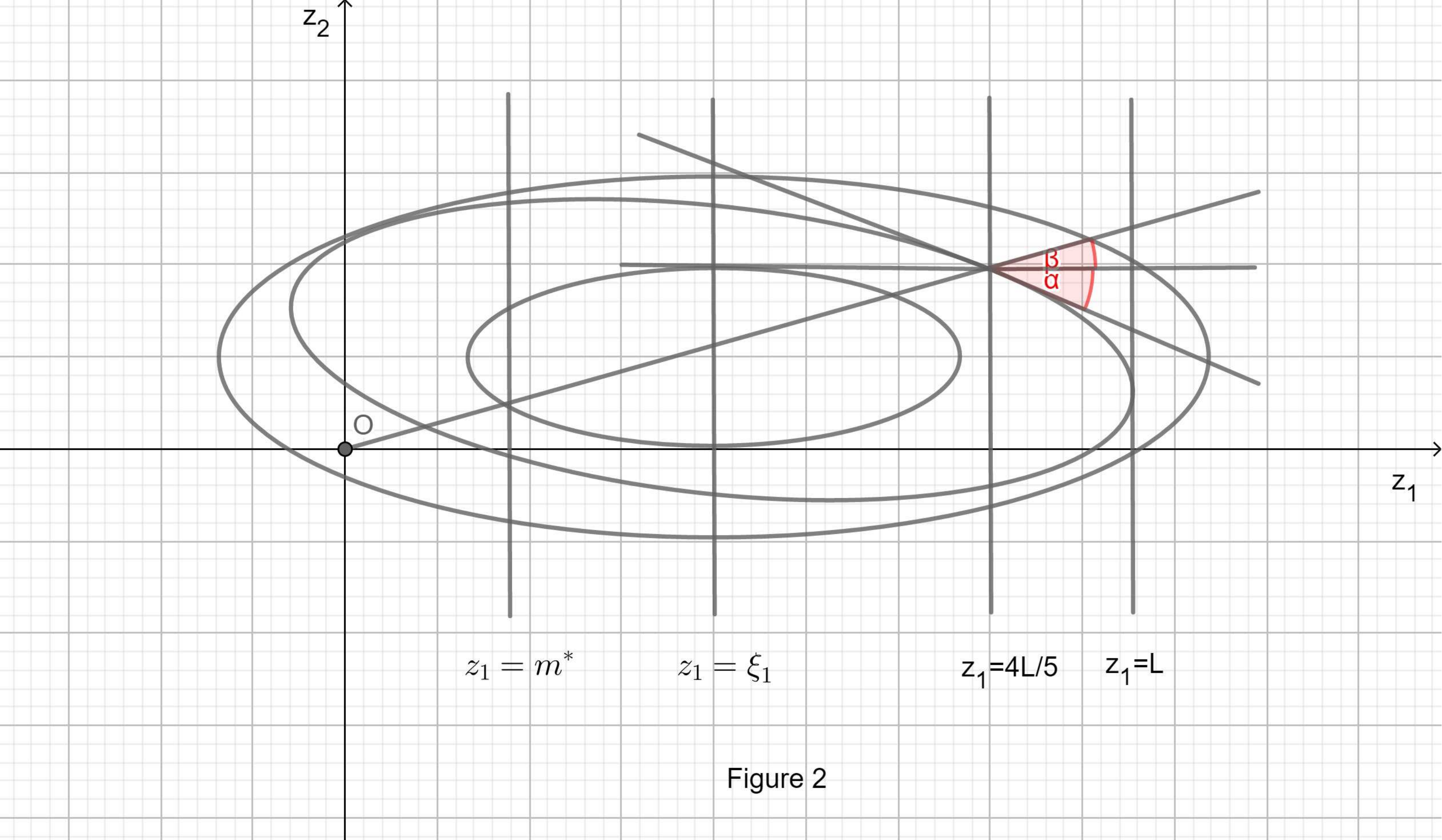}
  \caption{Estimation on lower bound of $l$}
\end{figure}

\vspace{40pt}

\section{Arbitrarily smallness of $u_{min}$ for $p\in(0,2)$}

Once we obtain the good shape estimation in Theorem \ref{t2.1}, it would be interesting to ask whether there holds an \textit{a-priori} lower bound for $u_{min}$. In this section, we shall show that this may not be true in some occasions for $p\in(0,2)$. We have the following theorem.

\begin{theo}\label{t3.1}
  For each $p\in(0,2)$, there exists a sequence of $f_j$ with uniformly upper and positive lower bound, which is also uniformly bounded in $C^\alpha({\mathbb{S}}^1)$ for some $\alpha\in(0,1)$, such that their corresponding positive classical solutions $u_j, j\in{\mathbb{N}}$ satisfy
    \begin{equation}\label{e3.1}
      \min_{{\mathbb{S}}^1}u_j\to0^+, \ \ \mbox{ as } j\to\infty.
    \end{equation}
\end{theo}

\noindent\textbf{Proof.} For each $\varepsilon>0$ small, we set
   $$
     u(\theta)=(\theta+\varepsilon)^{\frac{2}{2-p}}-\frac{2}{2-p}\varepsilon^{\frac{2}{2-p}-1}\theta, \ \ \theta\in[0,1]
   $$
to be an increasing positive function on $[0,1]$. Direct computation shows that
   \begin{eqnarray}\nonumber\label{e3.2}
     f(\theta)&\equiv&\frac{u_{\theta\theta}+u}{u^{p-1}}\\
     &=&\frac{a_1(\theta+\varepsilon)^{\frac{2(p-1)}{2-p}}+g(\theta)}{g^{p-1}(\theta)},
   \end{eqnarray}
where
   \begin{eqnarray*}
     a_1&\equiv&\frac{2}{2-p}\Big(\frac{2}{2-p}-1\Big),\\
     g(\theta)&\equiv&(\theta+\varepsilon)^{\frac{2}{2-p}}-\frac{2}{2-p}\varepsilon^{\frac{2}{2-p}-1}\theta.
   \end{eqnarray*}
Using the fact
   \begin{equation}\label{e3.3}
     \sigma_p(\theta+\varepsilon)^{\frac{2}{2-p}}\leq g(\theta)\leq C_p(\theta+\varepsilon)^{\frac{2}{2-p}}
   \end{equation}
for $p\in(0,2)$, where $\sigma_p$ is a constant closing to 1 and $C_p$ is large enough, one has
   \begin{equation}\label{e3.4}
     \frac{a_1+\sigma_p(\theta+\varepsilon)^2}{C_p^{p-1}}\leq f(\theta)\leq \frac{a_1+C_p(\theta+\varepsilon)^2}{\sigma_p^{p-1}}, \ \ \forall\theta\in[0,1]
   \end{equation}
and
   \begin{eqnarray}\nonumber\label{e3.5}
     u(1)&=&(1+\varepsilon)^{\frac{2}{2-p}}-\frac{2}{2-p}\varepsilon^{\frac{2}{2-p}-1}\\
     u_\theta(1)&=&\frac{2}{2-p}(1+\varepsilon)^{\frac{2}{2-p}-1}-\frac{2}{2-p}\varepsilon^{\frac{2}{2-p}-1}\\ \nonumber
     u_{\theta\theta}(1)&=&a_1(1+\varepsilon)^{\frac{2}{2-p}-2}.
   \end{eqnarray}
Next, let's construct a $C^{2,\alpha}$ function $\varphi$ on $[1,\pi]$ as follows. At first, we set
   \begin{eqnarray}\nonumber\label{e3.6}
     \varphi(1)&=&(1+\varepsilon)^{\frac{2}{2-p}}-\frac{2}{2-p}\varepsilon^{\frac{2}{2-p}-1}\\
     \varphi_\theta(1)&=&\frac{2}{2-p}(1+\varepsilon)^{\frac{2}{2-p}-1}-\frac{2}{2-p}\varepsilon^{\frac{2}{2-p}-1}\\ \nonumber
     \varphi_{\theta\theta}(1)&=&a_1(1+\varepsilon)^{\frac{2}{2-p}-2}.
   \end{eqnarray}
and let $\varphi_{\theta\theta}$ decreases rapidly to zero such that $\varphi_\theta$ approaches to $\frac{2}{2-p}$ on $(1,2.1)$. As a result, one has
   \begin{equation}\label{e3.7}
     \begin{cases}
       \varphi(2.1)\sim\frac{2.2}{2-p}+1, \ \ \varphi_\theta(2.1)=\frac{2}{2-p},\\
       \varphi_{\theta\theta}(2.1)=0.
     \end{cases}
   \end{equation}
Next, we let $\varphi_{\theta\theta}$ decreases rapidly to $-\frac{2}{2-p}$ such that $\varphi_\theta$ decreases to zero exactly at $\theta=\pi$. Furthermore, we assume that there holds
   \begin{equation}\label{e3.8}
     -\frac{2}{2-p}\leq\varphi_{\theta\theta}\leq0, \ \ \theta\in[2.1,\pi],
   \end{equation}
and thus implies that
   \begin{equation}\label{e3.9}
    \begin{cases}
     0\leq \varphi_\theta\leq\frac{2}{2-p}, & \forall\theta\in[2.1,\pi]\\
     a_2\leq\varphi\leq a_3, & \forall\theta\in[2.1,\pi]
    \end{cases}
   \end{equation}
for
   \begin{eqnarray*}
     a_2\equiv\frac{2.2}{2-p}, \ \ \ \ a_3\equiv\frac{2}{2-p}(\pi-1)+2.
   \end{eqnarray*}
Therefore, if one sets
   $$
    u(\theta)\equiv\begin{cases}
      g(\theta), & \theta\in[0,1)\\
       \varphi(\theta), & \theta\in(1,\pi]\\
       u(-\theta), & \theta\in(-\pi,0),
    \end{cases}
   $$
it is easy to verify that $u\in C^{2,\alpha}({\mathbb{S}}^1)$ for some $\alpha\in(0,1)$. Moreover, when $\theta\in(0,1]$, there holds \eqref{e3.4}. And when $\theta\in(1, 2.1]$,
   \begin{equation}\label{e3.10}
         g^{2-p}(1)\leq f(\theta)\leq\frac{\frac{2.2}{2-p}+1+a_1(1+\varepsilon)^{\frac{2}{2-p}-2}}{g^{p-1}(1)}.
   \end{equation}
Finally, if $\theta\in(2.1,\pi]$, one has that
   \begin{equation}\label{e3.11}
     \frac{\frac{0.2}{2-p}+1}{a_3^{p-1}} \leq f(\theta)\leq\frac{a_3}{a_2^{p-1}}
   \end{equation}
by \eqref{e3.8} and \eqref{e3.9}. The conclusion of the theorem follows by setting $\varepsilon=\frac{1}{j}$ and then letting $j\to+\infty$. $\Box$\\

As a corollary of Theorem \ref{t3.1}, one obtains Theorem \ref{t1.2}.

\vspace{40pt}

\section{Invertible Harnack inequality when $p\leq0$}

At the beginning of this section, we will first prove the following invertible Harnack inequality for non-positive $p$.

\begin{theo}\label{t4.1}
  Consider \eqref{e1.3} for $p\leq0$ and positive function $f\in C^\alpha({\mathbb{S}}^1)$. There exists a positive constant $C_{p,f}$ depending only on $p$ and $f_\pm$, such that
    \begin{equation}\label{e4.1}
     \begin{cases}
       C_{p,f}^{-1}u_{min}^p\leq u^2_{max}\leq C_{p,f} u_{min}^p, & \mbox{ for } p<0,\\
       C_{p,f}^{-1}\ln u_{min}^{-1}\leq u^2_{max}\leq C_{p,f}\ln u_{\min}^{-1}, & \mbox{ for } p=0
     \end{cases}
    \end{equation}
  holds for any positive classical solution $u$.
\end{theo}

Our proofs are consisted of two lemmas.

\begin{lemm}\label{l4.1}
  Under assumptions of Theorem \ref{t4.1}, there exists a positive constant $C_{p,f}$ such that
    \begin{equation}\label{e4.2}
     \begin{cases}
       u_{max}^2\leq C_{p,f}u_{min}^p, & \mbox{ for } p<0,\\
       u_{max}^2\leq C_{p,f}\ln u_{min}^{-1}, & \mbox{ for } p=0.
     \end{cases}
    \end{equation}
\end{lemm}

\noindent\textbf{Proof.} The lemma is a direct consequence of \eqref{e2.4} and \eqref{e2.8} for $p<0$, and a consequence of \eqref{e2.5} and \eqref{e2.8} for $p=0$. $\Box$\\

We have a second lemma under below, which can actually imply a stronger version of Theorem \ref{t4.1}.

\begin{lemm}\label{l4.2}
  Under the assumptions of Theorem \ref{t4.1}, there exists a positive constant $C_{p,f}$ such that
    \begin{equation}\label{e4.3}
     \begin{cases}
       u_{min}^{\frac{p}{2}}\leq C_{p,f}\int_{{\mathbb{S}}^1}u^{p-1}d\theta, & \mbox{ for } p<0,\\
       \sqrt{\ln u_{min}^{-1}}\leq C_{p,f}\int_{{\mathbb{S}}^1}u^{-1}d\theta, & \mbox{ for } p=0
     \end{cases}
    \end{equation}
  holds for any positive classical solution $u$ of \eqref{e1.3}.
\end{lemm}

\noindent\textbf{Proof.} Without loss of generality, one may assume that $u_{min}=u(0)$. By \eqref{e1.3}, there holds
   $$
    u_{\theta\theta}(\theta)\leq C_1u_{min}^{p-1}, \ \ \forall\theta\in[0,2\pi)
   $$
for some positive constant $C_1$. Thus, it follows from Taylor's expansion formula that
   \begin{equation}\label{e4.4}
      u(\theta)\leq u_{min}+C_1u_{min}^{p-1}\theta^2, \ \ \forall \theta\in[0,2\pi).
   \end{equation}
As a result, one obtains that
   \begin{eqnarray*}
     \int_{{\mathbb{S}}^1}u^{p-1}d\theta&\geq&\int_{{\mathbb{S}}^1}\Big(u_{min}+C_1u_{min}^{p-1}\theta^2\Big)^{p-1}d\theta\\
      &\geq& C_2\int_{{\mathbb{S}}^1}\Big(u_{min}^{1/2}+C_3u_{min}^{(p-1)/2}\theta\Big)^{2(p-1)}d\theta\\
      &\geq&\frac{C_2}{(2p-1)C_3}u_{min}^{(1-p)/2}\Big(u_{min}^{1/2}+C_3u_{min}^{(p-1)/2}\theta\Big)^{2p-1}\Big|_0^{u_{min}^{(1-p)/2}}\\
      &\geq& C_4u_{min}^{\frac{1-p}{2}}\Big(-C_{p,f}^{-1}+u_{\min}^{\frac{2p-1}{2}}\Big)\geq C_5u_{min}^{\frac{p}{2}}
   \end{eqnarray*}
when $p<0$, where \eqref{e4.4} and Lemma \ref{l2.2} have been used. When $p=0$, if $u_{min}$ is not small, then \eqref{e4.3} is clear true due to \eqref{e2.8}. Hence, one may assume that $u_{min}\ll1$. Note first that the function
   $$
    F(u)\equiv u^2-2f_+\ln u
   $$
is a decreasing function on $(0,\sqrt{f_+}]$ and is a increasing function on $[\sqrt{f_+},+\infty)$. Moreover,
  \begin{equation}\label{e4.5}
    \lim_{u\to 0^+}F(u)=\lim_{u\to+\infty}F(u)=+\infty.
  \end{equation}
If one denotes $\theta_2\in(0,2\pi)$ to be the first critical time of $u$, there must be
  \begin{equation}\label{e4.6}
     u(\theta_2)\geq\sqrt{f_-}
  \end{equation}
(since $u_{\theta\theta}$ is positive unless $u$ reaches $\sqrt{f_-}$ by \eqref{e1.3}) and
   \begin{equation}\label{e4.7}
     \theta_2\geq C_{p,f}^{-1}\frac{1}{\sqrt{\ln u_{min}^{-1}}},
   \end{equation}
where the second inequality follows from \eqref{e4.6} and the gradient bound
   \begin{equation}\label{e4.8}
     0<u_\theta\leq C_{p,f}\sqrt{\ln u_{min}^{-1}}, \ \ \forall \theta\in[0,\theta_2)
   \end{equation}
by \eqref{e2.4}. So, one has
    \begin{equation}\label{e4.9}
      u(\theta)\leq u_{min}+C_5\sqrt{\ln u_{min}^{-1}}\theta, \ \ \forall \theta\in[0,\theta_2)
    \end{equation}
and can conclude that
   \begin{eqnarray*}
     \int_{{\mathbb{S}}^1}u^{-1}d\theta&\geq&\int_{{\mathbb{S}}^1}\Big(u_{min}+C_5\sqrt{\ln u_{min}^{-1}}\theta\Big)^{-1}d\theta\\
      &\geq&\frac{1}{C_5\sqrt{\ln u_{min}^{-1}}}\ln\Big(u_{min}+C_5\sqrt{\ln u_{min}^{-1}}\theta\Big)\Big|^{C_{p,f}^{-1}\Big/\Big(2\sqrt{\ln u_{min}^{-1}}\Big)}_0\\
      &\geq&C_6\sqrt{\ln u_{min}^{-1}}.
   \end{eqnarray*}
The proof of the lemma was done. $\Box$\\

\noindent\textbf{Proof of Theorem \ref{t4.1}.} Integrating \eqref{e1.3} over ${\mathbb{S}}^1$, one gets that
  \begin{equation}\label{e4.10}
    \int_{{\mathbb{S}}^1}fu^{p-1}d\theta=\int_{{\mathbb{S}}^1}ud\theta.
  \end{equation}
Therefore, it yields from \eqref{e4.1}, \eqref{e4.3} and \eqref{e4.10} that
   \begin{equation}\label{e4.11}
       u_{min}^{\frac{p}{2}}\leq C_{p,f}\int_{{\mathbb{S}}^1}u^{p-1}d\theta\leq C_{p,f}\int_{{\mathbb{S}}^1} ud\theta\leq C_{p,f} u_{max}\leq C_{p,f}u_{min}^{\frac{p}{2}}
   \end{equation}
for $p<0$ and
   \begin{equation}\label{e4.12}
       \sqrt{\ln u_{min}^{-1}}\leq C_{p,f}\int_{{\mathbb{S}}^1}u^{-1}d\theta\leq C_{p,f}\int_{{\mathbb{S}}^1} ud\theta\leq C_{p,f} u_{max}\leq C_{p,f}\sqrt{\ln u_{min}^{-1}}
   \end{equation}
for $p=0$. The proof was done. $\Box$\\

As in Section 2, we want to prove the following round shape result.

\begin{theo}\label{t4.2}
  For each $p\in(-1,0]$ and positive function $f\in C^\alpha({\mathbb{S}}^1)$, there exists a positive constant $C_{p,f}$ depending only on $p$ and $f_\pm$, such that
    \begin{equation}\label{e4.13}
       C_{p,f}^{-1}\leq w_{\Omega_u}^-\leq w_{\Omega_u}^+\leq C_{p,f}w_{\Omega_u}^-
    \end{equation}
  holds for positive classical solution $u$ of \eqref{e1.3}.
\end{theo}

Adapting the notations and tricks as in Section 2, we will show the following proposition first.

\begin{prop}\label{p4.1}
  Under assumptions of Theorem \ref{t4.2}, there exists a positive constant $C_{p,f}\geq100$ such that
    \begin{equation}\label{e4.14}
      w_{\Omega_u}^+/w_{\Omega_u}^-\leq C_{p,f}.
    \end{equation}
\end{prop}

The next lemma will be used in the proofs of Proposition \ref{p4.1} and Theorem \ref{t4.2}.

\begin{lemm}\label{l4.3}
   There exists a positive constant $C_{p,f}$ such that
     \begin{equation}\label{e4.15}
        (w_{\Omega_u}^+)^{1-p}w_{\Omega_u}^-\geq C_{p,f}^{-1}
     \end{equation}
   holds for solution $u$ of \eqref{e1.3}.
\end{lemm}

\noindent\textbf{Proof.} Multiplying \eqref{e1.3} by $u$ and integrating over ${\mathbb{S}}^1$, the area $V(\Omega_u)$ of $\Omega_u$ satisfies that
   \begin{equation}\label{e4.16}
      V(\Omega_u)=\int_{{\mathbb{S}}^1}u(u_{\theta\theta}+u)d\theta=\int_{{\mathbb{S}}^1}fu^p\geq 2\pi f_-u_{max}^p\geq C_{p,f}^{-1}(w_{\Omega_u}^+)^p.
   \end{equation}
On the other hand,
   \begin{equation}\label{e4.17}
     V(\Omega_u)\leq C ab\leq Cw_{\Omega_u}^+w_{\Omega_u}^-
   \end{equation}
holds for $a,b$ defined in \eqref{e2.10} and universal constant $C>0$. Hence, \eqref{e4.14} follows from \eqref{e4.15} and \eqref{e4.16}. $\Box$\\

\noindent\textbf{Proof of Proposition \ref{p4.1}.} Suppose on the contrary, then
  \begin{equation}\label{e4.18}
    w_{\Omega_u}^+/w_{\Omega_u}^-\geq C_{p,f}\geq100.
  \end{equation}
When $m^*\leq3L/4$, if $p<0$, there hold
  \begin{equation}\label{e4.19}
    l^2(\theta_1)-l^2(\pi/2)\leq\frac{2f_+}{p}\Big(u^p(\theta_1)-u^p(\pi/2)\Big)\leq Cl^p
  \end{equation}
and \eqref{e2.20}. So, we have
   \begin{equation}\label{e4.20}
      L^2\leq C_{p,f}l^p.
   \end{equation}
Combining \eqref{e4.20} with \eqref{e4.15} and \eqref{e4.18}, we conclude that
   \begin{equation}\label{e4.21}
      C_{p,f}^{-1}\leq l\leq L\leq C_{p,f}
   \end{equation}
and hence \eqref{e4.14} thanks to $p\in(-1,0)$. If $p=0$, we use
  $$
   l^2(\theta_1)-l^2(\pi/2)\leq2f_+\ln\frac{u(\theta_1)}{u(\pi/2)}\leq0
  $$
and \eqref{e2.20} to conclude contradiction.

If $m^*\geq 3L/4$, we define $\theta_3\in(\pi/2,\pi)$ by ${r}(\theta_3)=\Big(L/2, \varphi(L/2)\Big)$ and define $\theta_4\in(\pi,2\pi)$ by ${r}(\theta_4)=(3L/4,\varphi(3L/4))$. Moreover, we denote
  $$
    \delta_1\equiv u(\theta_3), \ \ \delta_2\equiv u(\theta_4).
  $$
Similarly, one needs also consider the lower portion of boundary $\partial\Omega$ and denote its graph by
   $$
     \widehat{\varphi}(s):\ {r(\theta)}=(s,\widehat{\varphi}(s)), \ \ s\in[\pi,2\pi).
   $$
 we define $\widehat{\theta}_3\in(\pi/2,\pi)$ by ${r}(\widehat{\theta}_3)=\Big(L/2, \widehat{\varphi}(L/2)\Big)$ and define $\widehat{\theta}_4\in(\pi,2\pi)$ by ${r}(\widehat{\theta}_4)=(3L/4,\widehat{\varphi}(3L/4))$. Parallelly, we can also define
  $$
    \widehat{\delta}_1\equiv u(\widehat{\theta}_3), \ \ \widehat{\delta}_2\equiv u(\widehat{\theta}_4)
  $$
and $m_*$ to be the $z_1$ coordinate of ${r}(3\pi/2)$.

  The remaining proofs are divided into several lemmas.

\begin{lemm}\label{l4.4}
   Under the assumptions of Proposition \ref{p4.1} and \eqref{e4.18}, there exists a small constant $\varepsilon\equiv\varepsilon_{p,f}\in(0,1)$ such that if $m^*\geq 3L/4$, there holds
     \begin{equation}\label{e4.22}
        \delta_1\equiv u(\theta_3)\leq\begin{cases}
           \varepsilon L^{\frac{2}{p}}, &  p<0,\\
           le^{-\varepsilon L^2}, & p=0.
        \end{cases}
     \end{equation}
   And if $m^*\leq 3L/4$, there holds
      \begin{equation}\label{e4.23}
         \delta_2\equiv u(\theta_4)\leq\begin{cases}
           \varepsilon L^{\frac{2}{p}}, &  p<0,\\
           le^{-\varepsilon L^2}, & p=0.
        \end{cases}
      \end{equation}
   Similarly, a same result holds for $m^*$ is replaced by $m_*$ and $\delta_1, \delta_2$ is replaced by $\widehat{\delta}_1,\widehat{\delta}_2$.
\end{lemm}

\noindent\textbf{Proof.} When $m^*\geq 3L/4$, arguing as in Section 2 and using the geometric lemma \ref{l2.3}, if \eqref{e4.22} is not true, one has
  \begin{equation}\label{e4.24}
    l^2(\pi/2)-l^2(\theta_3)\leq\frac{2f_+}{p}\Big(u^p(\pi/2)-u^p(\theta_3)\Big)\leq \varepsilon L^2
  \end{equation}
in case $p<0$ and
   \begin{equation}\label{e4.25}
     l^2(\pi/2)-l^2(\theta_3)\geq(3L/4)^2-\Big((L/2)^2+l^2\Big)\geq L^2/5
   \end{equation}
by assumption \eqref{e4.18}. Contradiction holds. If $p=0$, we use
   $$
    l^2(\pi/2)-l^2(\theta_3)\leq 2f_+\ln\frac{u(\pi/2)}{u(\theta_3)}\leq0
   $$
and \eqref{e4.25} to achieve a similar contradiction.

When $m^*\leq 3L/4$, one can use
  \begin{equation}\label{e4.26}
    l^2(\theta_1)-l^2(\theta_4)\leq\frac{2f_+}{p}\Big(u^p(\theta_1)-u^p(\theta_4)\Big)\leq C_{p,f}u^p(\theta_4)
  \end{equation}
and
   \begin{equation}\label{e4.27}
    l^2(\theta_1)-l^2(\theta_4)\geq C_{p,f}^{-1}L^2
   \end{equation}
to conclude the validity of \eqref{e4.23} in case $p<0$. For $p=0$,
  $$
   l^2(\theta_1)-l^2(\theta_4)\leq 2f_+\ln\frac{u(\theta_1)}{u^p(\theta_4)}
  $$
is used to replace \eqref{e4.26}. Together with \eqref{e4.27}, a contradiction yields. The inequalities for lower portion can be proved similarly. $\Box$\\

Considering an arc interval
   $$
    \omega\equiv\Big\{\theta\in(\pi/2,\pi]|\ {r}(\theta)=(s,\varphi(s)), -d\leq s\leq L/2\Big\}=(\theta_3,\pi)
   $$
and defining
   $$
      \beta(\theta)\equiv\arg({r}(\theta))\in{\mathbb{S}}^1,\ \  \theta\in{\mathbb{S}}^1
   $$
to be the reverse Gaussian mapping $G^{-1}: {\mathbb{S}}^1\to{\mathbb{S}}^1$, one has
   \begin{equation}\label{e4.28}
     \tan(\theta_3-\pi/2)\leq \tan\beta(\theta_3)\leq\frac{l}{L/2}\ll1
   \end{equation}
and hence
    \begin{equation}\label{e4.29}
       \Big\{\theta|\ 51\pi/100\leq\theta\leq\pi\Big\}\subset\omega.
    \end{equation}
Using the expansion relation \eqref{e2.1}, direct computation shows that
   \begin{equation}\label{e4.30}
      \frac{d\beta}{d\theta}=\frac{u(u+u_{\theta\theta})}{l^2}, \ \ l^2=u^2+u_\theta^2.
   \end{equation}
Therefore, it is inferred from \eqref{e1.3} and \eqref{e4.29} that
   \begin{equation}\label{e4.31}
      \int_{G^{-1}(\omega)}u^{-p}(\beta)l^2(\beta)d\beta=\int_{\omega}fd\theta\geq 49\pi f_-/100
   \end{equation}
where we denote
    $$
     u(\beta)\equiv u(G(\beta)), \ \ r(\beta)\equiv r(G(\beta))
    $$
for short. Now, we can deduce a contradiction by estimating L.H.S. of \eqref{e4.31} using a similar geometric lemma as in \cite{CL}.

\begin{lemm}\label{l4.5}
  Under the assumptions of Proposition \ref{p4.1}, \eqref{e4.18} and $m^*\geq 3L/4$, there exists a positive constant $C_{p,f}$ such that
    \begin{equation}\label{e4.32}
       \int_{G^{-1}(\omega)}u^{-p}(\beta)r^{2}(\beta)d\beta\leq\begin{cases} C_{p,f}\Big(L^{\frac{2}{p}+1-p}+L^{\frac{2(1-p)}{p}}\Big), & p<0\\
        C_{p,f}\Big(L^{2-p}e^{-\varepsilon L^2}+(Le^{-\varepsilon L})^{2-p}\Big), & p=0.
       \end{cases}
    \end{equation}
\end{lemm}

\noindent\textbf{Proof.} We will follow the arguments in \cite{CL} to prove \eqref{e4.32} Denoting $l_P$ to be the tangential line of $\partial\Omega$ at $P\equiv{r}(\theta_3)$, and denoting $l_P^\perp=\overrightarrow{OQ}$ to be the perpendicular ray of $l_P$ starting from origin which intersect with $l_P$ at $Q$, we set
   $$
    \omega_1\equiv\Big\{\beta\in G^{-1}(\omega)|\ {r}(G(\beta))\in\triangle OPQ\Big\}
   $$
and
   $$
    \omega_2\equiv\Big\{\beta\in G^{-1}(\omega)|\ \cos\beta\in[-d,0]\Big\}.
   $$
it is clear that
   \begin{equation}\label{e4.33}
     G^{-1}(\omega)\subset\omega_1\cup\omega_2
   \end{equation}
and
   \begin{equation}\label{e4.34}
       \cup_{t\in[0,1], \beta\in\omega_1}\Big\{t{r}(\beta)\Big\}\subset\triangle OPQ, \ \
       \cup_{t\in[0,1], \beta\in\omega_2}\Big\{t{r}(\beta)\Big\}\subset{\mathcal{R}},
   \end{equation}
where ${\mathcal{R}}$ is the rectangular
   $$
    {\mathcal{R}}\equiv\Big\{(z_1,z_2)\in{\mathbb{R}}^2|\ z_1\in[-d,0], z_2\in[-\widehat{h},h]\Big\},
   $$
where
    $$
     h\equiv\varphi(0), \ \ \widehat{h}\equiv-\widehat{\varphi}(0).
    $$
\begin{figure}
  \includegraphics[width=5.1in,height=2.55in]{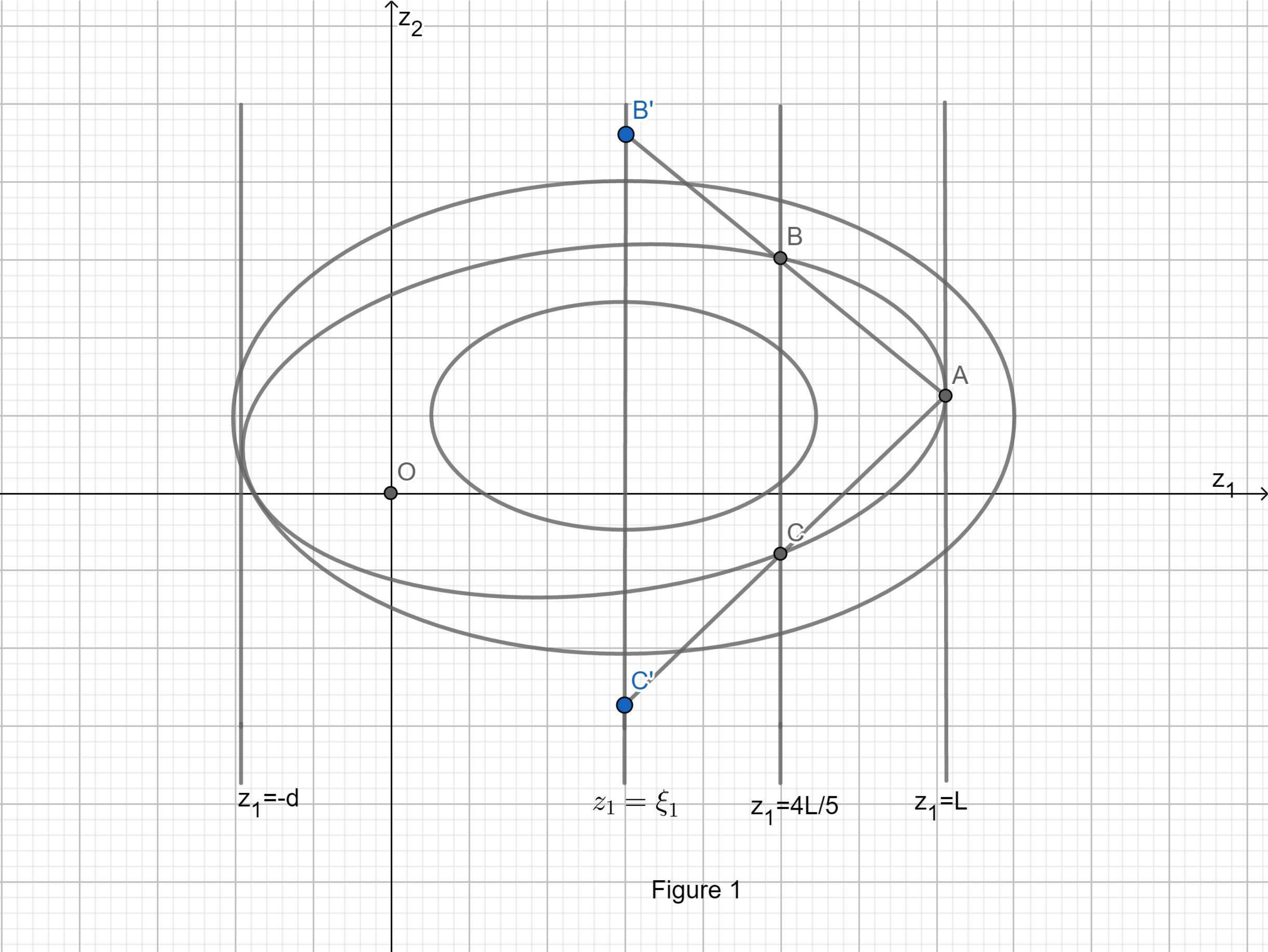}
  \caption{Estimation on lower bound of $l$}
\end{figure}
For any $\beta\in\omega_1$, let's define
   $$
    \widetilde{r}\equiv\sup\Big\{tr(\beta)|\ t{r}(\beta)\in\triangle OPQ\Big\}.
   $$
For $\beta\in\omega_2$, we can also define
   $$
    \widetilde{r}\equiv\sup\Big\{tr(\beta)|\ t{r}(\beta)\in{\mathcal{R}}\Big\}.
   $$
Now, we go to estimate
  \begin{equation}\label{e4.35}
    \int_{\omega_1}\widetilde{r}^{2-p}(\beta)d\beta+\int_{\omega_2}\widetilde{r}^{2-p}(\beta)d\beta\geq\int_{G^{-1}(\omega)}u^{-p}r^{2}(\beta)d\beta.
  \end{equation}
At first, noting that by \eqref{e4.22}, there holds
   \begin{equation}\label{e4.36}
     h=\varphi(0)\leq 2|OQ|=2\delta_1\leq \begin{cases}
       2\varepsilon L^{\frac{2}{p}}, & \forall p<0\\
       2le^{-\varepsilon L^2}, & \forall p=0.
     \end{cases}
   \end{equation}
Similarly, one also has
   \begin{equation}\label{e4.37}
     \widehat{h}=-\widehat{\varphi}(0)\leq\begin{cases}
       2\varepsilon L^{\frac{2}{p}}, & \forall p<0\\
       2le^{-\varepsilon L^2}, & \forall p=0.
     \end{cases}
   \end{equation}
On another hand, by the convexity of $\Omega$, the triangle $\triangle_0ABC$ composed of
  $$
   A=(-d,\varphi(-d)), \ \ B=\xi-be_2, \ \ C=\xi+be^2
  $$
lies entirely inside of $\Omega$. So, one has
  $$
   \frac{d}{L/2}\leq\frac{h+\widehat{h}}{2b}
  $$
and thus
  \begin{equation}\label{e4.38}
    d\leq\begin{cases}
      \varepsilon L^{1+\frac{2}{p}}l^{-1}\leq\varepsilon L^{1+\frac{2}{p}}u_{min}^{-1}\leq\varepsilon L^{1+\frac{2}{p}}u_{max}^{-\frac{2}{p}}\leq\varepsilon L, & p<0\\
      CLe^{-\varepsilon L^2}, & p=0
    \end{cases}
  \end{equation}
by invertible Harnack inequality \eqref{e4.1}. As a result, we could estimate
  \begin{eqnarray*}
    \int_{\omega_1}\widetilde{r}^{2-p}(\beta)d\beta&\leq&\int^{\frac{\pi}{2}-\arcsin\frac{\delta_1}{\sqrt{(L/2)^2-l^2}}}_0\Bigg(\frac{\delta_1}{\cos\beta}\Bigg)^{2-p}d\beta\\
     &\leq& C\delta_1^{2-p}\int^{\frac{\pi}{2}-\arcsin\frac{\delta_1}{\sqrt{(L/2)^2+l^2}}}_0\Big(\frac{\pi}{2}-\beta\Big)^{p-2}d\beta\\
     &\leq& C_{p,f}\delta_1L^{1-p}\leq\begin{cases} C_{p,f}L^{\frac{2}{p}+1-p}, & p<0\\
       C_{p,f}L^{1-p}le^{-\varepsilon L^2}, & p=0
     \end{cases}
  \end{eqnarray*}
and
  \begin{eqnarray*}
    \int_{\omega_2}\widetilde{r}^{2-p}(\beta)d\beta&\leq&4 \int^{\arcsin\frac{d}{\sqrt{d^2+(h+\widehat{h})^2}}}_0\Big(\frac{h}{\cos\beta}\Big)^{2-p}d\beta\\
     &\leq& C_{p,f}(h^{2-p}+hd^{1-p})\leq\begin{cases}
       C_{p,f}\Big(L^{\frac{2(2-p)}{p}}+L^{\frac{2}{p}+1-p}\Big), & p<0\\
       C_{p,f}(Le^{-\varepsilon L})^{2-p}, & p=0.
     \end{cases}
  \end{eqnarray*}
The proof Lemma \ref{l4.5} was completed. $\Box$\\

\noindent\textbf{Continue the proof of Proposition \ref{p4.1}.} By Lemma \ref{l4.5} and noting $p\in(-1,0]$, it follows from \eqref{e4.31} and \eqref{e4.32} that
    $$
     C_{p,f}^{-1}\leq\int_{G^{-1}(\omega)}u^{-p}r^{2}(\beta)d\beta\leq C_{p,f}L^{-\sigma_p}
    $$
for some $\sigma_p>0$. So, there holds
   $$
    C_{p,f}^{-1}\leq l\leq C_{p,f}L\leq C_{p,f}
   $$
by \eqref{e4.15}. The conclusion of Proposition \ref{p4.1} was drawn. $\Box$\\

\noindent\textbf{Complete the proof of Theorem \ref{t4.2}.} By Lemma \ref{l4.3} and Proposition \ref{p4.1}, one concludes that
  \begin{equation}\label{e4.39}
    C_{p,f}^{-1}\leq l\leq C_{p,f} L\leq C_{p,f}l
  \end{equation}
for some positive constant $C_{p,f}$. So, Theorem \ref{t4.2} follows from \eqref{e2.10}, \eqref{e2.11} and \eqref{e4.39}. $\Box$\\

\vspace{40pt}

\section{Large body and droplet argument}

In the previous section we obtain a round shape result. In order to deduce a positive, two-sided bound on the support function, we need to exclude the possibility of large bodies. Here we show

\begin{theo}\label{t5.1}
  Under the assumptions of Theorem \ref{t4.2}, there exists a positive constant $C_{p,f}$ depending only on $p, \min f$ and $||f||_{C^\alpha({\mathbb{S}}^n)}$, such that
     \begin{equation}\label{e5.1}
        w_{\Omega_u}^+\leq C_{p,f}
     \end{equation}
  holds for positive classical solution $u$ of \eqref{e1.3}. As a result, there holds
     \begin{equation}\label{e5.2}
       C_{p,f}^{-1}\leq u_{min}\leq u_{max}\leq C_{p,f}.
     \end{equation}
\end{theo}

Inequality \eqref{e5.2} is a direct consequence of \eqref{e5.1} by invertible Harnack inequality \eqref{e4.1}. Therefore, the crucial step in proof of Theorem \ref{t5.1} consists of establishing \eqref{e5.1}. Under below, We will first blow down the solution to a limiting convex body with droplet shape and then produce a contradiction.\\

\noindent\textbf{Proof.} Supposing that \eqref{e5.1} is not true for some fixed $p\in(-1,0]$, then given any $k\in{\mathbb{N}}$, there exist a positive function $f_k\in C^{\alpha}({\mathbb{S}}^1), \alpha\in(0,1)$ satisfying
   \begin{equation}\label{e5.3}
      C_0^{-1}\leq f_k\leq C_0,\ \ ||f_k||_{C^\alpha({\mathbb{S}}^1)}\leq C_0,  \ \ \forall k\in{\mathbb{N}}
   \end{equation}
for some positive $C_0$ independent of $k$, and a positive classical solution $u_k$ to \eqref{e1.3} such that
   \begin{equation}\label{e5.4}
      \max_{{\mathbb{S}}^1}u_k=u_k(0)\equiv\lambda_k\uparrow+\infty, \ \ \forall k,
   \end{equation}
by rotation if necessary. Applying Theorem \ref{t4.2} to $u_k$, one has
   \begin{equation}\label{e5.5}
     C_{p,C_0}^{-1}\leq w_{\Omega_{u_k}}^-\leq w_{\Omega_{u_k}}^+\leq C_{p,C_0}w_{\Omega_{u_k}}^-
   \end{equation}
for some positive constant $C_{p,C_0}$ independent of $k$. Rescaling $u_k$ by the function $v_k=\lambda_k^{-1}u_k$, we get a solution to
   \begin{equation}\label{e5.6}
     \partial_\theta^2v_k+v_k=\lambda_k^{p-2}f_kv_k^{p-1}
   \end{equation}
which satisfies that
   \begin{equation}\label{e5.7}
     v_k(0)=\max_{{\mathbb{S}}^1}v_k=1, \ \ \forall k.
   \end{equation}
Another hand, since $\Omega_{v_k}=\lambda_k^{-1}\Omega_{u_k}$, there holds
   \begin{equation}\label{e5.8}
     w_{\Omega_{v_k}}^+=\lambda_k^{-1}w_{\Omega_{u_k}}^+, \ \ w_{\Omega_{v_k}}^-=\lambda_k^{-1}w_{\Omega_{u_k}}^-.
   \end{equation}
Noting that
    \begin{equation}\label{e5.9}
      w_{\Omega_{u_k}}/2\leq\lambda_k\leq w_{\Omega_{u_k}},
    \end{equation}
it follows from \eqref{e5.5} and \eqref{e5.9} that
   \begin{equation}\label{e5.10}
      C_1^{-1}\leq w_{\Omega_{v_k}}^-\leq w_{\Omega_{v_k}}^+\leq C_1
   \end{equation}
for some positive constant $C_1$ independent of $k$. Multiplying \eqref{e5.6} by $v_k$, integrating over ${\mathbb{S}}^1$ and then performing integration by parts, one gets that
  \begin{equation}\label{e5.11}
   \int_{{\mathbb{S}}^1}|\partial_\theta v_k|^2=\int_{{\mathbb{S}}^1}v_k^2-\lambda_k^{p-2}\int_{{\mathbb{S}}^1}f_kv_k^p\leq 2\pi, \ \ \forall k.
  \end{equation}
With the help of Sobolev embedding theorem, we conclude that
   \begin{equation}\label{e5.12}
     ||v_k||_{C^{1/2}({\mathbb{S}}^1)}\leq C_2, \ \ \forall k
   \end{equation}
for some positive constant $C_2$ independent of $k$. By Arzela-Ascoli theorem, there exists a limiting nonnegative function $v_\infty\in C^{1/2}({\mathbb{S}}^1)$, such that for a subsequence of $k$,
   \begin{equation}\label{e5.13}
    \begin{cases}
     v_k\to v_\infty, \ \ \mbox{ uniformly on }  {\mathbb{S}}^1,\\
     v_k\to v_\infty, \ \ \mbox{ uniformly on }  C^{2,\alpha}(P),\ \ P\equiv\Big\{\theta\in[0,2\pi)|\ v_\infty(\theta)>0\Big\},\\
     v_\infty(0)=\max_{{\mathbb{S}}^1}v_\infty=1, \ \ \min_{{\mathbb{S}}^1}v_\infty=0,\ \ v_\infty\in H^1({\mathbb{S}}^1),\\
     \partial_\theta^2v_\infty+v_\infty=0, \ \  \forall \theta\in P.
    \end{cases}
   \end{equation}
Noting that the unique solution to
   $$
    \partial^2_\theta v_\infty+v_\infty=0, \ \ v_\infty(0)=\max_{{\mathbb{S}}^1}v_\infty=1
   $$
is given by
   \begin{equation}\label{e5.14}
     v_\infty(\theta)=\cos\theta, \ \ \forall\theta\in(-\pi/2,\pi/2).
   \end{equation}
We claim that
   \begin{equation}\label{e5.15}
     v_\infty(\theta)\equiv-\kappa\cos\theta, \ \ \forall\theta\in[\pi/2,3\pi/2]
   \end{equation}
for some nonnegative constant $\kappa$. In fact, noting that
    $$
      \int_{{\mathbb{S}}^1}|\partial_\theta u_k|^2\leq\int_{{\mathbb{S}}^1} u_k^2\Rightarrow \int_{{\mathbb{S}}^1}|\partial_\theta v_k|^2\leq\int_{{\mathbb{S}}^1} v_k^2,
    $$
by lower semi-continuity of weak convergence on $H^1({\mathbb{S}}^1)$, one concludes that
   \begin{equation}\label{e5.16}
     \int_{{\mathbb{S}}^1}|\partial_\theta v_\infty|^2\leq\int_{{\mathbb{S}}^1} v_\infty^2.
   \end{equation}
Using \eqref{e5.14}, there holds
   \begin{equation}\label{e5.17}
     \int^{3\pi/2}_{\pi/2}|\partial_\theta v_\infty|^2\leq \int^{3\pi/2}_{\pi/2}v_\infty^2.
   \end{equation}
However, since
    $$
      v_\infty(\pi/2)=v_\infty(3\pi/2)=0,
    $$
it is inferred from \eqref{e5.17} and Wirtinger's inequality that
   $$
     v_\infty(\theta)=-\kappa\cos\theta, \ \ \forall\theta\in[\pi/2,3\pi/2].
   $$
The claim was done. On another hand, sine the convex bodies $\Omega_{v_k}$ sub-converges to the convex body $\Omega_{v_\infty}$ in weak sense, one has also
  \begin{equation}\label{e5.18}
     \Omega_{v_\infty}\supseteq E_\infty\equiv\Bigg\{z=(z_1,z_2)\in{\mathbb{R}}^2\Big|\ \frac{(z_1-\xi_{1\infty})^2}{a_\infty^2}+\frac{(z_2-\xi_{2\infty})^2}{b_\infty^2}=1\Bigg\}
  \end{equation}
by \eqref{e5.10}, where $\xi_\infty=(\xi_{1\infty},\xi_{2\infty})\in{\mathbb{R}}^2$ and $a_\infty\geq b_\infty$ are positive constants. However, a convex body $\Omega_{v_\infty}$ with support function
   $$
    v_\infty(\theta)=\begin{cases}
      \cos\theta, & \forall\theta\in(-\pi/2,\pi/2),\\
      -\kappa\cos\theta, & \forall\theta\in[\pi/2,3\pi/2]
    \end{cases}
   $$
can only be a degenerate two-sides thin droplet with zero in-radius, which contradiction with \eqref{e5.18}. The proof of Theorem \ref{t5.1} was completed. $\Box$\\

\vspace{40pt}

\section{Topological degree and Theorem \ref{t1.3}}

In this section, we shall use the \textit{a-priori} upper/lower bound of positive classical solution of \eqref{e1.3} to prove the desired solvability result Theorem \ref{t1.3}. The method of Leray-Schauder's topological degree has been used by Chou-Wang in \cite{CW} for $L_p$-Minkowski problem and later developed by
Chen-Li in \cite{CL} to dual-Minkowski problem. For the convenience of the reader, we will present a proof here. For each $p\in(-1,0]$ and positive function $f\in C^\alpha({\mathbb{S}}^1)$, let's denote
   $$
   p_t\equiv tp, \ \ f_t\equiv tf+(1-t), \ \ \forall t\in[0,1]
   $$
and consider the equation
   \begin{equation}\label{e6.1}
     u_{\theta\theta}+u=f_tu^{p_t-1}, \ \ \forall\theta\in{\mathbb{S}}^1.
   \end{equation}
Setting
   \begin{equation}\label{e6.2}
    \Gamma\equiv\Big\{t\in[0,1]|\ \mbox{there is a positive classical solution to } \eqref{e6.1}\Big\},
   \end{equation}
 we turn to prove $\Gamma=[0,1]$. As shown in Section 8, the uniqueness of \eqref{e1.3} with $f\equiv1$ was not true for the case $p\in(-1,0)$. Thus, to show the desired solvability for $-1<p\leq0$, we need to utilize a uniqueness result by Chow \cite{C} for $n\geq1, p=1-n$ and $f\equiv1$.

\begin{lemm}\label{l6.1}
  For each $n\geq1$ and $p=1-n$, the solution of \eqref{e1.1} with respect to constant function $f\equiv1$ is unique.
\end{lemm}

Now, let's apply the method of topological degree to prove Theorem \ref{t1.3} in case of $-1<p\leq0$. As above, for each $p\in(-1,0]$ and positive function $f\in C^\alpha({\mathbb{S}}^1)$, we consider the equation \eqref{e6.1}.

\begin{prop}\label{p6.1}
  For each $-1<p\leq0$ and $f$ satisfying assumptions of Theorem \ref{t5.1}, there holds $\Gamma=[0,1]$ for $\Gamma$ defined by \eqref{e6.2}.
\end{prop}

\noindent\textbf{Proof.} Let's use topological degree to show $\Gamma=[0,1]$. At first, we define
   $$
     {\mathcal{F}}_t(u)\equiv u_{\theta\theta}+u-f_tu^{p_t-1}, \ \ \forall t\in[0,1].
   $$
By Theorem \ref{t5.1} and Schauder's estimates for linear elliptic partial differential equations, there exists a positive constant $C_*$ independent of $t$, such that
   \begin{equation}\label{e6.3}
     C_*^{-1}\leq u\leq C_*, \ \ ||u||_{C^{2,\alpha}({\mathbb{S}}^1)}\leq C_*, \ \ \alpha\in(0,1)
   \end{equation}
holds for each zero $u$ of ${\mathcal{F}}_t$. So, if one defines $$
 {\mathcal{O}}\equiv\Big\{u\in C^{2,\alpha}({\mathcal{S}}^1)|\  (2C_*)^{-1}\leq u\leq 2C_*, \ \ ||u||_{C^{2,\alpha}({\mathbb{S}}^1)}\leq C_*\Big\},
$$
it is clear that
   $$
    {\mathcal{F}}_t^{-1}(0)\cap\partial{\mathcal{O}}=0.
   $$
Noting that by Chow's uniqueness result Lemma \ref{l6.1}, $u_0\equiv1$ is the unique solution to \eqref{e6.1} with respect to $f_0\equiv1$. Moreover, the linearized equation
   \begin{equation}\label{e6.4}
     \varphi_{\theta\theta}+\varphi=-\varphi, \ \ \forall\theta\in{\mathbb{S}}^1
   \end{equation}
of \eqref{e6.1} at $t=0, u_0\equiv1$ has only trivial solution $\varphi\equiv0$ since the unique $2\pi$-periodic function
  $$
   \varphi(\theta)=A\cos\sqrt{2}\theta+B\sin\sqrt{2}\theta
  $$
is given by $A=B=0$. Thus, the degree
   $$
    deg({\mathcal{F}}_0,{\mathcal{O}},0)=1\not=0.
   $$
Combining with the preservation property of topological degree \cite{L,N2}, there holds
   $$
    deg({\mathcal{F}}_t,{\mathcal{O}},0)=deg({\mathcal{F}}_0,{\mathcal{O}},0)\not=0, \ \ \forall t\in[0,1].
   $$
So, the solvability of \eqref{e6.1} for $t\in[0,1]$ holds true. The proof of Proposition \ref{p6.1} was done. $\Box$\\

\vspace{40pt}

\section{Trigonometric identity for $p\leq-2$}

In this section, we turn to prove Theorem \ref{t1.4}. First, we establish a crucial trigonometric identity for \eqref{e1.3}.

\begin{lemm}\label{l7.1}
  Given $p\leq-2$ and a nonnegative function $f\in C^\alpha({\mathbb{S}}^1)$ which is piece wise $C^1$, then for any positive classical solution $u$ of \eqref{e1.3}, there holds
    \begin{equation}\label{e7.1}
       \int_{{\mathbb{S}}^1}K_{f}(\theta)u^p=0,
    \end{equation}
  where
     \begin{equation}\label{e7.2}
        K_f(\theta)\equiv(p+2)f\cos2\theta+f_\theta\sin2\theta.
     \end{equation}
\end{lemm}

\noindent\textbf{Proof.} Multiplying \eqref{e1.3} by $\cos2\theta u$, integrating over ${\mathbb{S}}^1$ and then forming integration by parts, one gets that
  \begin{equation}\label{e7.3}
   -\int_{{\mathbb{S}}^1}\cos2\theta u_\theta^2-\int_{{\mathbb{S}}^1}\cos2\theta u^2=\int_{{\mathbb{S}}^1}\cos2\theta fu^p.
  \end{equation}
Multiplying again \eqref{e1.3} by $\sin2\theta u_\theta$, one gets that
  \begin{equation}\label{e7.4}
    -\int_{{\mathbb{S}}^1}\cos2\theta u_\theta^2-\int_{{\mathbb{S}}^1}\cos2\theta u^2=-\frac{1}{p}\int_{{\mathbb{S}}^1}\Big[f(\theta)\sin2\theta\Big]_\theta u^p.
  \end{equation}
Thus, there holds
   $$
     \int_{{\mathbb{S}}^1}\Bigg\{pf(\theta)\cos2\theta+\Big[f(\theta)\sin2\theta\Big]_\theta\Bigg\}u^p=0,
   $$
which is equivalent to \eqref{e7.1} and \eqref{e7.2}. The proof was done. $\Box$\\

In this section, we always denote
   $$
    t^\beta\equiv|t|^{\beta-1}t, \ \ \forall t\not=0
   $$
and use the relations like
   $$
    t^{\beta-1}\cdot t^{-\beta}=|t|^{-1}\not=t^{-1}
   $$
to distinguish the usually ones.\\

\noindent\textbf{Proof of Theorem \ref{t1.4}.} We consider first the case $p=-2$.  Letting $f(\theta)=2+\cos2\theta$ and supposing there is a positive classical solution $u$ of \eqref{e1.3}, one obtains that
   $$
    K_f(\theta)=-2\sin^22\theta\leq0, \ \ K_f\not\equiv0, \ \  \forall\theta\in{\mathbb{S}}^1
   $$
by Lemma \ref{l7.1}. So, it yields from \eqref{e7.1} a contradiction
   $$
     0>\int_{{\mathbb{S}}^1}K_f(\theta)u^p=0.
   $$
The proof for $p=-2$ was completed. To show the non-existence result for $p<-2$, let's first introduce a positive function by
   \begin{equation}\label{e7.5}
    \xi(\theta)\equiv\begin{cases}
      \int^{\pi/4}_\theta \sin^{\frac{p}{2}}2\vartheta d\vartheta, & \theta\in(0,\pi/2)\\
      \xi(-\theta), & \theta\in(-\pi/2,0)\\
      \xi(\theta-\pi), & \theta\in(\pi/2,3\pi/2).
    \end{cases}
   \end{equation}
It's easy to see that $\xi$ is a $\pi-$periodic even function on
  $$
     {\mathbb{R}}\setminus\{k\pi/2, k\in{\mathbb{Z}}\Big\}.
  $$
Moreover, for $p<-2$, there holds
   $$
    \lim_{\theta\to\frac{k\pi}{2}}\xi(\theta)=+\infty
   $$
for each $k\in{\mathbb{Z}}$. Now, we need a second lemma.

\begin{lemm}\label{l7.2}
  For each $p<-2$, the function
    $$
     \phi(\theta)\equiv\begin{cases}
       |\sin2\theta|^{-\frac{p+2}{2}}\xi(\theta), & \forall\theta\not=\frac{k\pi}{2}, k\in{\mathbb{Z}},\\
       -\frac{1}{p+2}, & \forall\theta=k\pi, k\in{\mathbb{Z}},\\
       \frac{1}{p+2}, & \forall\theta=\frac{\pi}{2}+k\pi, k\in{\mathbb{Z}}
     \end{cases}
    $$
  is a $\pi-$periodic even $Lip({\mathbb{S}}^1)$ function satisfying
     \begin{equation}\label{e7.7}
       \phi(\theta)>\frac{1}{p+2}, \ \ \forall\theta\in{\mathcal{P}}\equiv[0,2\pi)\setminus\{\pi/2+k\pi\}.
     \end{equation}
\end{lemm}

\noindent\noindent\textbf{Proof.} Since $\phi$ is a $\pi-$periodic even function and Lipschitz continuity everywhere $\theta\not=k\pi/2, k\in{\mathbb{Z}}$, let us first show $\phi$ is differentiable at $\theta=0$. The case $\theta=\pi/2$ is similarly. The proof is elementary by L'Hospital's law. Let's first prove $\phi$ is continuous at $\theta=0$. In fact,
  \begin{eqnarray*}
    \lim_{\theta\to0}\phi(\theta)&=&\lim_{\theta\to0}\frac{\int^{\pi/4}_\theta\sin^{\frac{p}{2}}2\vartheta d\vartheta}{(2\theta)^{\frac{p+2}{2}}}\\
     &=&\lim_{\theta\to0}-\frac{\sin^{\frac{p}{2}}2\theta}{(p+2)(2\theta)^{\frac{p}{2}}}=-\frac{1}{p+2}.
  \end{eqnarray*}
So $\phi$ is continuous. Next, we show that $\phi$ is differentiable at $\theta=0$. In fact, since
   \begin{eqnarray*}
     \lim_{\theta\to0}\frac{\phi(\theta)+\frac{1}{p+2}}{\theta}&=&\lim_{\theta\to0}\frac{\int^{\pi/4}_{\theta}\sin^{\frac{p}{2}}2\vartheta d\vartheta+\frac{1}{p+2}\sin^{\frac{p+2}{2}}2\theta}{\theta\sin^{\frac{p+2}{2}}2\theta}\\
     &=&\lim_{\theta\to0}\frac{-\sin^{\frac{p}{2}}2\theta+\sin^{\frac{p}{2}}2\theta\cos2\theta}{(p+4)2^{\frac{p}{2}}\theta^{\frac{p+2}{2}}}\\
     &=&\lim_{\theta\to0}\frac{-1+\cos2\theta}{(p+4)\theta}=0,
   \end{eqnarray*}
$\phi$ is a $Lip({\mathbb{S}}^1)$ function which is greater than $\frac{1}{p+2}$ everywhere $\theta\in{\mathcal{P}}$. $\Box$\\

Now, we define
   $$
     f(\theta)\equiv-\frac{1}{p+2}+\phi(\theta), \ \ \forall\theta\in{\mathbb{S}}^1
   $$
to be a positive Lipschitz function outside two polar of ${\mathbb{S}}^1$. Direct calculation shows that
  \begin{eqnarray}\nonumber\label{e7.8}
    K_f(\theta)&=&(p+2)\Bigg(-\frac{1}{p+2}+\sin^{-\frac{p+2}{2}}2\theta\int^{\pi/4}_{\theta}\sin^{\frac{p}{2}}2\vartheta d\vartheta\Bigg)\cos2\theta\\
    &&-(p+2)\sin^{-\frac{p+4}{2}}2\theta\cos2\theta\int^{\pi/4}_{\theta}\sin^{\frac{p}{2}}2\vartheta d\vartheta\cdot\sin2\theta-1\\ \nonumber
    &=&-1-\cos2\theta<0, \ \ \forall\theta\in(0,\pi/2).
  \end{eqnarray}
Since $K_f$ is $\pi$-periodical even function, \eqref{e7.8} is in conflict with \eqref{e7.1}. The proof of Theorem \ref{t1.4} was completed. $\Box$\\

\vspace{40pt}

\section{Uniqueness for constant $f$}

In this section, let's complete the proof of Theorem \ref{t1.5}. We assume that $u(0)=u_{min}$ by rotation if necessary and denote it to be $m$ for short. Multiplying \eqref{e1.3} by $u_\theta$, integrating over ${\mathbb{S}}^1$ and then performing integration by parts, it yields that
   \begin{equation}\label{e8.1}
     u_\theta^2+u^2-\frac{2}{p}u^p\equiv m^2-\frac{2}{p}m^p, \ \ \forall\theta\in{\mathbb{S}}^1
   \end{equation}
 for $p\not=0$ and
   \begin{equation}\label{e8.2}
    u_\theta^2+u^2-2\ln u\equiv m^2-2\ln m, \ \ \forall\theta\in{\mathbb{S}}^1
   \end{equation}
 for $p=0$. Setting
    $$
     F(u)\equiv \begin{cases}
        u^2-\frac{2}{p}u^p, & p\not=0,\\
        u^2-2\ln u, & p=0,
     \end{cases}
    $$
 it is clear that $F$ is monotone decreasing on $(0,1]$ and monotone increasing on $[1,+\infty)$. Moreover,
   \begin{eqnarray}\label{e8.3}
     &&\lim_{u\to0+} F(u)=\begin{cases}
        0, & p\in(0,2)\\
        +\infty, & p\leq0,
     \end{cases}\\ \nonumber
    &&\lim_{u\to+\infty} F(u)=+\infty, \ \ \forall p<2, p\not=0.
   \end{eqnarray}
 On the other hand, it follows from the maximum principle that for non-constant solution $u$, there holds
   \begin{equation}\label{e8.4}
     m=\min_{{\mathbb{S}}^1}u<1, \ \ M\equiv\max_{{\mathbb{S}}^1}u>1.
   \end{equation}
 By \eqref{e8.1}-\eqref{e8.2} and picture of $F$, it is clear that
    \begin{equation}\label{e8.5}
      F(m)=F(M), \ \ u\uparrow \mbox{ from } m \mbox{ to } M, \mbox{ and } u\downarrow \mbox{ from } M \mbox{ to }m.
    \end{equation}
 Furthermore, using the uniqueness of first order ordinary differential equation \eqref{e8.1}-\eqref{e8.2} and a reflection $v(\theta)=u(-\theta)$, one can deduce that $u$ must be symmetric around $\theta=\pi$ and thus
     \begin{equation}\label{e8.6}
       u_\theta(\pi)=0.
     \end{equation}
 Therefore, it is inferred from \eqref{e8.1}-\eqref{e8.2} that any possible non-constant solution of \eqref{e1.3} must satisfy the compatible condition
   \begin{equation}\label{e8.7}
     \begin{cases}
       \int^M_m\frac{du}{\sqrt{F(m)-F(u)}}=\pi/\kappa,\\
       F(m)=F(M)
     \end{cases}
   \end{equation}
for some positive integer $\kappa$. In fact, one has the following proposition.

 \begin{prop}\label{p8.1}
   Supposing that $p<2$ and $f\equiv1$, \eqref{e1.3} has a positive classical solution $u$ satisfying $u_{min}=m\in(0,1)$ if and only if there exists a pair $(m,M), M>1$ satisfying \eqref{e8.7} with some $\kappa\in{\mathbb{N}}$.
 \end{prop}

For any $m\in(0,1)$ and letting $M=M(m)>1$ being the unique positive constant determined by the second relation in \eqref{e8.7}, we define
   $$
    H(m)\equiv\int^{M(m)}_m\frac{du}{\sqrt{F(m)-F(u)}}.
   $$
As shown below, $H$ is a $C^1$-function on $(0,1)$. Furthermore, the following lemma would be useful in proving of Theorem \ref{t1.5}.

\begin{lemm}\label{l8.1} Letting $p<2$, there holds
   \begin{eqnarray*}
     \lim_{m\to0^+}H(m)&=&\begin{cases}
        \frac{\pi}{2-p}, & \forall p\in(0,2)\\
        \pi/2, & \forall p\in(-\infty,0],
     \end{cases}\\
     \lim_{m\to1^-}H(m)&=&\frac{\pi}{\sqrt{2-p}}, \ \ \forall p<2.
   \end{eqnarray*}

\end{lemm}

\noindent\textbf{Proof.} We consider the case $0<p<2$ first. Since $M(m)\to\Big(\frac{2}{p}\Big)^{1/(2-p)}$ as $m$ tends to zero, one has
   \begin{eqnarray*}
     \lim_{m\to0^+}H(m)&=&\int^{\big(\frac{2}{p}\big)^{1/(2-p)}}_0\frac{du}{\sqrt{\frac{2}{p}u^p-u^2}}\\
      &=&\frac{2}{2-p}\arcsin\frac{u^{\frac{2-p}{2}}}{\sqrt{2/p}}\Bigg|^{\big(\frac{2}{p}\big)^{1/(2-p)}}_0=\frac{\pi}{2-p}.
   \end{eqnarray*}
To calculate the limit at $m=1$, we first use the Cauchy's mean value theorem to deduce that
  $$
    \frac{(m^2-2m)-(u^2-2u)}{(m^2-\frac{2}{p}m^p)-(u^2-\frac{2}{p}u^p)}=\frac{2\xi-2}{2\xi-2\xi^{p-1}}\sim\frac{1}{2-p}, \ \ \xi\in(m,u).
  $$
Therefore,
  \begin{eqnarray*}
    \lim_{m\to1^-}H(m)&=&\lim_{m\to1^-}\frac{1}{\sqrt{2-p}}\int^{M(m)}_m\frac{du}{\sqrt{(m-1)^2-(u-1)^2}}\\
     &=&\lim_{m\to1^-}-\frac{1}{\sqrt{2-p}}\arcsin\frac{u-1}{m-1}\Big|^{M(m)}_m\\
     &=&\lim_{m\to1^-}\frac{1}{\sqrt{2-p}}\Bigg(\arcsin\frac{M(m)-1}{1-m}+\frac{\pi}{2}\Bigg)=\frac{\pi}{\sqrt{2-p}}
  \end{eqnarray*}
for any $p<2$, where
    \begin{eqnarray*}
      &&\lim_{M\to1^+}\frac{(M-1)^2}{M^2-\frac{2}{p}M^p-F(1)}=\lim_{m\to1^-}\frac{(m-1)^2}{m^2-\frac{2}{p}m^p-F(1)}=\frac{1}{2-p}, \ \ p\not=0\\
      &&\lim_{M\to 1^+}\frac{(M-1)^2}{M^2-2\ln M-1}=\lim_{m\to1^-}\frac{(m-1)^2}{m^2-2\ln m-1}=\frac{1}{2}, \ \ p=0
    \end{eqnarray*}
and second relation of \eqref{e8.7} have been used. Now, let's calculate the limit at $m=0$ for $p\leq0$. In fact, for each fixed $\varepsilon>0$, if $m$ is chosen small, $\frac{1}{\sqrt{F(m)-F(u)}}$ is also small on the fixed interval $[0,\varepsilon^{-1}]$. On another hand, one also has
   \begin{eqnarray*}
    F(m)&=&\begin{cases}
       (1+o_\varepsilon(1))\Big(-\frac{2}{p}m^p-u^2\Big), & p<0\\
       (1+o_\varepsilon(1))\Big(-2\ln m-u^2\Big), & p=0
    \end{cases}\\
    &=&(1+o_\varepsilon(1))(M^2-u^2), \ \ \forall u\geq\varepsilon^{-1},
   \end{eqnarray*}
where $o_\varepsilon(1)$ is a small quantity as long as $\varepsilon$ is small. So, we obtain that
   \begin{eqnarray*}
    H(m)&=&o_m(1)+(1+o_\varepsilon(1))\int^{M(m)}_{\varepsilon^{-1}}\frac{du}{\sqrt{M^2-u^2}}\\
    &&+\begin{cases}
      (1+o_\varepsilon(1))\sqrt{\frac{|p|}{2}}\int^\varepsilon_m\frac{du}{\sqrt{m^p-u^p}}, & p<0\\
            (1+o_\varepsilon(1))\sqrt{\frac{1}{2}}\int^\varepsilon_m\frac{du}{\sqrt{\ln m^{-1}-\ln u^{-1}}}, & p=0
      \end{cases}\\
    &=&o_m(1)+(1+o_\varepsilon(1))\Bigg(\arcsin\frac{u}{M(m)}\Big|^{M(m)}_{\varepsilon^{-1}}+R(m,\varepsilon)\Bigg)\\
    &=&o_{m,\varepsilon}(1)+(1+o_\varepsilon(1))\Big(\pi/2-\arcsin\frac{\varepsilon^{-1}}{M(m)}\Big)
   \end{eqnarray*}
for small quantity $o_{m,\varepsilon}(1)$ with respect to $m, \varepsilon$, where
  \begin{eqnarray*}
    R(m,\varepsilon)&\equiv&\sqrt{\frac{|p|}{2}}\int^\varepsilon_m\frac{du}{\sqrt{m^p-u^p}}\leq\sqrt{\frac{|p|}{2}}\varepsilon^{\frac{2-p}{2}}\int^\varepsilon_m\frac{u^{\frac{p-2}{2}}du}{\sqrt{m^p-u^p}}\\
     &=&-\sqrt{\frac{2}{|p|}}\varepsilon^{\frac{2-p}{2}}\arcsin{\Big(\frac{u}{m}\Big)^{\frac{p}{2}}}\Big|^\varepsilon_m=\sqrt{\frac{2}{|p|}}\varepsilon^{\frac{2-p}{2}}\Bigg(\frac{\pi}{2}-\arcsin\Big(\frac{\varepsilon}{m}\Big)^{\frac{p}{2}}\Bigg)=o_\varepsilon(1)
  \end{eqnarray*}
for $p<0$ and
  \begin{eqnarray*}
    R(m,\varepsilon)&\equiv&\sqrt{\frac{1}{2}}\int^\varepsilon_m\frac{du}{\sqrt{\ln u-\ln m}}\leq\sqrt{\frac{1}{2}}\varepsilon\sqrt{\ln\varepsilon^{-1}}\int^\varepsilon_m\frac{\frac{1}{u\sqrt{-\ln u}}du}{\sqrt{\ln u-\ln m}}\\
     &=&-\varepsilon\sqrt{\ln\varepsilon^{-1}}\arcsin\sqrt{\frac{\ln u}{\ln m}}\Big|^{\varepsilon}_m=\varepsilon\sqrt{\ln{\varepsilon^{-1}}}\Bigg(\frac{\pi}{2}-\arcsin\sqrt{\frac{\ln\varepsilon}{\ln m}}\Bigg)=o_\varepsilon(1)
  \end{eqnarray*}
for $p=0$ have been used. Thus, the limit
   $$
    \lim_{m\to0^+}H(m)=\pi/2
   $$
was shown for $p\leq0$ and the conclusion was drawn. $\Box$\\

To proceed further, let's differentiate the second relation of \eqref{e8.7} on $m$ and derive that
    \begin{equation}\label{e8.8}
      \frac{dM}{dm}=\frac{m-m^{p-1}}{M-M^{p-1}}.
    \end{equation}
To differentiate the first relation on \eqref{e8.7}, let's first perform integration by parts to yields
   \begin{eqnarray*}
     H(m)&=&\int^M_{1+\sigma}\frac{du}{\sqrt{F(M)-F(u)}}+\int^{1-\sigma}_{m}\frac{du}{\sqrt{F(m)-F(u)}}+\int^{1+\sigma}_{1-\sigma}\frac{du}{\sqrt{F(m)-F(u)}}\\
      &=&\int^M_{1+\sigma}\sqrt{F(M)-F(u)}d\frac{1}{u-u^{p-1}}+\int^{1-\sigma}_m\sqrt{F(m)-F(u)}d\frac{1}{u-u^{p-1}}\\
      &&+\frac{\sqrt{F(m)-F(1+\sigma)}}{(1+\sigma)-(1+\sigma)^{p-1}}-\frac{\sqrt{F(m)-F(1-\sigma)}}{(1-\sigma)-(1-\sigma)^{p-1}}+\int^{1+\sigma}_{1-\sigma}\frac{du}{\sqrt{F(m)-F(u)}}.
   \end{eqnarray*}
 Differentiating on $m$, there holds
   \begin{eqnarray}\nonumber\label{e8.9}
     \frac{dH}{dm}&=&\int^M_{1+\sigma}\frac{m-m^{p-1}}{\sqrt{F(m)-F(u)}}d\frac{1}{u-u^{p-1}}+\int^{1-\sigma}_{m}\frac{m-m^{p-1}}{\sqrt{F(m)-F(u)}}d\frac{1}{u-u^{p-1}}\\
     &&+\frac{m-m^{p-1}}{(1+\sigma)-(1+\sigma)^{p-1}}\frac{1}{\sqrt{F(m)-F(1+\sigma)}}-\int^{1+\sigma}_{1-\sigma}\frac{(m-m^{p-1})du}{(F(m)-F(u))^{3/2}}\\ \nonumber
     &&-\frac{m-m^{p-1}}{(1-\sigma)-(1-\sigma)^{p-1}}\frac{1}{\sqrt{F(m)-F(1-\sigma)}}\\ \nonumber
     &&\equiv (m-m^{p-1})\Big(J_1+J_2+J_3+J_4+J_4\Big)
   \end{eqnarray}
for a fixed small constant $\sigma>0$. Setting $G(u)\equiv F(u)-F(1)$, it yields that
   \begin{eqnarray*}
    && \int^{1+\sigma}_1\frac{du}{(F(m)-F(u))^{\frac{3}{2}}}=\int^{1+\sigma}_1\frac{du}{(G(m)-G(u))^{\frac{3}{2}}}=\int^{1+\sigma}_1\frac{\frac{1}{G^{\frac{3}{2}}(u)}}{\Big(\frac{G(m)}{G(u)}-1\Big)^{\frac{3}{2}}}du\\
      && \ \ \ \ \ \ \ \ \ \ \ \ \ \ =-\frac{1}{G(m)}\int^{1+\sigma}_1\frac{\frac{\sqrt{G(u)}}{G'(u)}d\frac{G(m)}{G(u)}}{\Big(\frac{G(m)}{G(u)}-1\Big)^{\frac{3}{2}}}=\frac{2}{G(m)}\int^{1+\sigma}_1\frac{\sqrt{G(u)}}{G'(u)}d\frac{1}{\sqrt{\frac{G(m)}{G(u)}-1}}\\
      &=&\frac{2G(1+\sigma)}{G(m)G'(1+\sigma)}\frac{1}{\sqrt{G(m)-G(1+\sigma)}}-\frac{2}{G(m)}\int^{1+\sigma}_1\frac{\frac{1}{2}-\frac{G(u)}{(G'(u))^2}G''(u)}{\sqrt{G(m)-G(u)}}du,
   \end{eqnarray*}
where
    $$
     \lim_{u\to1^+}\frac{\sqrt{G(u)}}{G'(u)}=\frac{1}{2\sqrt{2-p}}
    $$
has been used. Similarly, there holds
   \begin{eqnarray*}
     \int^1_{1-\sigma}\frac{du}{(F(m)-F(u))^{\frac{3}{2}}}&=&-\frac{2G(1-\sigma)}{G(m)G'(1-\sigma)}\frac{1}{\sqrt{G(m)-G(1-\sigma)}}\\
     &&-\frac{2}{G(m)}\int^1_{1-\sigma}\frac{\frac{1}{2}-\frac{G(u)}{(G'(u))^2}G''(u)}{\sqrt{G(m)-G(u)}}du.
   \end{eqnarray*}
So,
  \begin{eqnarray}\nonumber\label{e8.10}
   J_3+J_4+J_5&=&\frac{2\sqrt{G(m)-G(1+\sigma)}}{G'(1+\sigma)G(m)}-\frac{2\sqrt{G(m)-G(1-\sigma)}}{G'(1-\sigma)G(m)}\\
   &&+\frac{2}{G(m)}\int^{1+\sigma}_{1-\sigma}\frac{\frac{1}{2}-\frac{G(u)}{(G'(u))^2}G''(u)}{\sqrt{G(m)-G(u)}}du.
  \end{eqnarray}
On another hand, there hold
  \begin{eqnarray*}
    J_1&=&-2\int^M_{1+\sigma}\frac{G''(u)/(G'(u))^2}{\sqrt{G(m)-G(u)}}du=\frac{2}{G(m)}\int^M_{1+\sigma}\frac{\frac{G^{3/2}(u)G''(u)}{(G'(u))^3}d\frac{G(m)}{G(u)}}{\sqrt{\frac{G(m)}{G(u)}-1}}\\
    &=&\frac{4}{G(m)}\int^M_{1+\sigma}\frac{G^{3/2}(u)G''(u)}{(G'(u))^3}d\sqrt{\frac{G(m)}{G(u)}-1}\\
    &=&-\frac{4}{G(m)}\frac{G(1+\sigma)G''(1+\sigma)}{(G'(1+\sigma))^3}\sqrt{G(m)-G(1+\sigma)}\\
    &&-\frac{4}{G(m)}\int^M_{1+\sigma}\frac{GG'G'''+\frac{3}{2}(G')^2G''-3G(G'')^2}{(G')^4}\sqrt{G(m)-G(u)}du
  \end{eqnarray*}
and
   \begin{eqnarray*}
    J_2&=&\frac{4}{G(m)}\frac{G(1-\sigma)G''(1-\sigma)}{(G'(1-\sigma))^3}\sqrt{G(m)-G(1-\sigma)}\\
    &&-\frac{4}{G(m)}\int^{1-\sigma}_m\frac{GG'G'''+\frac{3}{2}(G')^2G''-3G(G'')^2}{(G')^4}\sqrt{G(m)-G(u)}du
   \end{eqnarray*}
Summation as above, one concludes that
   \begin{eqnarray*}
     \frac{dH/dm}{m-m^{p-1}}&=&J_1+J_2+J_3+J_4+J_5\\
     &=&\frac{-4G(1+\sigma)G''(1+\sigma)+2(G'(1+\sigma))^2}{G(m)(G'(1+\sigma))^3}\sqrt{G(m)-G(1+\sigma)}\\
     &&+\frac{4G(1-\sigma)G''(1-\sigma)-2(G'(1-\sigma))^2}{G(m)(G'(1-\sigma))^3}\sqrt{G(m)-G(1-\sigma)}\\
     &&-\frac{4}{G(m)}\int_{[m,M]\setminus(1-\sigma,1+\sigma)}\frac{GG'G'''+\frac{3}{2}(G')^2G''-3G(G'')^2}{(G')^4}\sqrt{G(m)-G(u)}du\\
     &&+\frac{2}{G(m)}\int^{1+\sigma}_{1-\sigma}\frac{\frac{1}{2}-\frac{G(u)}{(G'(u))^2}G''(u)}{\sqrt{G(m)-G(u)}}du.
   \end{eqnarray*}
Using again integration by parts,
  \begin{eqnarray*}
    &&\frac{2}{G(m)}\int^{1+\sigma}_{1-\sigma}\frac{\frac{1}{2}-\frac{G(u)}{(G'(u))^2}G''(u)}{\sqrt{G(m)-G(u)}}du=\frac{-2(G')^2+4GG''}{G(m)(G')^3}\sqrt{G(m)-G(u)}\Big|^{1+\sigma}_{1-\sigma}\\
     && \ \ \ \ \ \ \ \ \ \ \ \ \ \ \ \ \ \ \ +\frac{2}{G(m)}\int^{1+\sigma}_{1-\sigma}\Bigg[\frac{(G')^2-2GG''}{(G')^3}\Bigg]_u\sqrt{G(m)-G(u)}du.
  \end{eqnarray*}
Hence, we arrive at the following relation formula.

\begin{prop}\label{p8.2}
 For each $m\in(0,1)$ and $M(m)$ determined by second relation of \eqref{e8.7}, the $C^1(0,1)-$function defined by
   $$
    H(m)\equiv\int^{M(m)}_m\frac{du}{\sqrt{F(m)-F(u)}}, \ \ \forall 0<m<1
   $$
satisfies an intrinsic relation
   \begin{eqnarray}\label{e8.11}
    -\frac{G(m)dH/dm}{4(m-m^{p-1})}=\int^M_m\frac{K(u)}{(G'(u))^4}\sqrt{G(m)-G(u)}du
   \end{eqnarray}
with the kernel
   $$
    K(u)\equiv GG'G'''+\frac{3}{2}(G')^2G''-3G(G'')^2.
   $$
\end{prop}

\vspace{10pt}

\noindent\textbf{Remark 8.1} It is remarkable to note that $u=1$ is a removable singularity of the function $K(u)/(G'(u))^4$ as shown below. Another hand, after integration by parts, one can obtain another intrinsic identity
  \begin{equation}\label{e8.12}
     \frac{G(m)dH/dm}{m-m^{p-1}}=\int^M_m\frac{(G')^2-2GG''}{(G')^2}\frac{1}{\sqrt{G(m)-G(u)}}du
  \end{equation}
using the relation
   \begin{equation}\label{e8.13}
    -2\frac{K(u)}{(G'(u))^4}=\Bigg[\frac{(G')^2-2GG''}{(G')^3}\Bigg]'_u.
   \end{equation}

\noindent To proceed further, we note that by direct calculation, the kernel
  $$
   K(u)\equiv GG'G'''+\frac{3}{2}(G')^2G''-3G(G'')^2
  $$
satisfies that
  \begin{eqnarray}\nonumber\label{e8.14}
    K'(u)&=&\frac{4(p-1)(p-2)}{pu^{1-p}}L(u),\\
    \Big(u^3L'(u)\Big)'&=&2(p-2)u^{p-3}T(u),\\ \nonumber
    T'(u)/p&=&2u\Big\{(p-1)(3p-4)u^{p-2}-(2p^2-7p+8)\Big\}
  \end{eqnarray}
where
   \begin{eqnarray*}
     L(u)&\equiv&(3p-4)u^{2p-4}+2(2p-1)(p-2)u^{p-4}\\
    &&+(p-2)(p-8)u^{-2}-2(2p^2-7p+8)u^{p-2}-p(p-3)\\
    T(u)&\equiv&2(p-1)(3p-4)u^{p}+(p-2)(p-4)(2p-1)\\
    &&-p(2p^2-7p+8)u^{2}
   \end{eqnarray*}
satisfies that
  \begin{equation}\label{e8.15}
    K(1)=L(1)=L'(1)=T(1)=0.
  \end{equation}
As a result, we have the following monotone result of $H$.
\begin{lemm}\label{l8.2}
  For $p\in(1,2)$, the function $H(m)$ is a monotone decreasing function on $m\in(0,1)$. While, for $p\in[1/2,1)$, the function $H(m)$ is a monotone increasing function on $m\in(0,1)$.
\end{lemm}

\noindent\textbf{Proof.} When $1\leq p\leq 4/3$, there holds
 \begin{eqnarray}\nonumber\label{e8.16}
   T'(u)/p&=&2u\Big\{(p-1)(3p-4)u^{p-2}-(2p^2-7p+8)\Big\}\\
   &\leq& u(p-2)(3p^2-5p+8)<0, \ \ \forall u\in\Bigg(0,\Bigg(\frac{2}{p}\Bigg)^{\frac{1}{2-p}}\Bigg)
 \end{eqnarray}
So, we conclude that $K/(G')^4$ is a continuous negative function  on
  $$
    \Bigg(0,\Bigg(\frac{2}{p}\Bigg)^{\frac{1}{2-p}}\Bigg)\setminus\{1\}
  $$
for $1<p\leq4/3$. Thus, the lemma follows by \eqref{e8.11}.

For $p\in(0,2)\setminus[1,4/3]$, noting that \eqref{e8.16} and \eqref{e8.15},
one has that
    \begin{equation}\label{e8.17}
     T'(u)/p\leq 2u(p^2-4)\leq 0, \ \ \forall u>1\Rightarrow T(u)<0, \ \ \forall u>1.
    \end{equation}
Another hand, since $T'$ vanishes at a unique point $u_0\in(0,1)$, and
    $$
     T'(u)=\begin{cases}
        >0, & \forall u\in(0,u_0)\\
        <0, & \forall u>u_0,
     \end{cases}
    $$
there holds
    \begin{eqnarray}\nonumber\label{e8.18}
      T(u)&\geq&\min\Big\{T(0),T(1)\Big\}\\
      &=&\min\Big\{(p-2)(p-4)(2p-1),0\Big\}=0, \ \ \forall u\in(0,1)
    \end{eqnarray}
as long as $p\in[1/2,2)$. Consequently, it follows from $K(1)=K'(1)=L(1)=L'(1)=0$ and \eqref{e8.17}-\eqref{e8.18} that
     \begin{equation}\label{e8.19}
       K(u)<0, \ \ \forall u\in\Big(0,(2/p)^{\frac{1}{2-p}}\Big)\setminus\{1\}
     \end{equation}
in case of $p\in(1,2)$ and
     \begin{equation}\label{e8.20}
       K(u)>0, \ \ \forall u\in\Big(0,(2/p)^{\frac{1}{2-p}}\Big)\setminus\{1\}
     \end{equation}
in case of $p\in[1/2,1)$. Hence, we obtain that the monotonicity of $H$ by \eqref{e8.11} and \eqref{e8.19}-\eqref{e8.20}. $\Box$\\

We also have the following monotonicity of $H$ at the end point $m=1$.

\begin{lemm}\label{l8.3}
   For $p\in(-\infty,2)\setminus\{1,-2\}$, one has
     \begin{equation}\label{e8.21}
       \frac{dH}{dm}=\begin{cases}
         <0, & p\in(1,2)\\
         >0, & p\in(-2,1)\\
         <0, & p\in(-\infty,-2)
       \end{cases}
     \end{equation}
for $0<1-m\ll1$.
\end{lemm}

\noindent\textbf{Proof.} To calculate the sign of derivative of $H$ near $m=1$, we use the first relation \eqref{e8.11}. In fact, for $m$ closing to $1$, one has
  \begin{eqnarray*}
     -\frac{G(m)dH}{4(m-m^{p-1})dm}&=&\int^M_m\frac{K(u)}{(G'(u))^4}\sqrt{G(m)-G(u)}du\\
      &=&\begin{cases}
         <0, & p\in(1,2)\\
         >0, & p\in(-2,1)\\
         <0, & p\in(-\infty,-2)
      \end{cases}
   \end{eqnarray*}
for all $0<1-m\ll1$, where
   \begin{eqnarray*}
    \lim_{u\to1}T'(u)/p=2u(p^2-4)\begin{cases}
        <0, & p\in(-2,2)\\
        >0, & p\in(-\infty,-2)
    \end{cases}
   \end{eqnarray*}
have been used. The proof of the lemma was done. $\Box$\\

\noindent\textbf{Remark 8.2} For any $\gamma\in(0,G(m))$, we denote $u\in(m,1), v\in(1,M(m))$ by
    $$
     G(u)=G(v)=\gamma,
    $$
where $M(m)$ is determined by second relation in \eqref{e8.7}. It is inferred from \eqref{e8.12} that
  \begin{equation}\label{e8.22}
    \frac{G(m)dH/dm}{m-m^{p-1}}=-\int^{G(m)}_{0}\mathcal{K}(u,v)\frac{d\gamma}{\sqrt{G(m)-\gamma}},
  \end{equation}
where
   $$
    \mathcal{K}(u,v)\equiv\frac{(G'(u))^2-2G(u)G''(u)}{(G'(u))^3}-\frac{(G'(v))^2-2G(v)G''(v)}{(G'(v))^3}.
   $$
Therefore, one gets the following proposition.

\begin{prop}\label{p8.3}
 Suppose that for some $p<2, p\not=1$,
   \begin{equation}\label{e8.23}
    {\mathcal{K}}(u,v)\ \mbox{ doesn't change sign on curve } {\mathcal{L}}
  \end{equation}
 holds for curve defined by
   \begin{equation}\label{e8.24}
      {\mathcal{L}}\equiv\Big\{(u,v)\in{\mathbb{R}}^2|\ G(u)=G(v), \ \ u\in(m,1), v\in(1,M(m))\Big\}.
   \end{equation}
We have $H$ is a monotone function on $m\in(0,1)$.
\end{prop}

It is remarkable that Proposition \ref{p8.3} reduces an uniqueness problem to a purely algebraic condition \eqref{e8.23} on $G(\cdot)$. Moreover, it is not hard to verify that for $p\in[1/2,2)$, the assumption of Proposition \ref{p8.3} holds true since
   $$
    \frac{(G')^2-2GG''}{(G')^3}
   $$
is a monotone function on $(m,M(m))$ by \eqref{e8.13} and the sign of $K$ in proof of Lemma \ref{l8.2}.\\

 Now, let's complete the proof of Theorem \ref{t1.5} as below.

\noindent\textbf{Proof of Theorem \ref{t1.5}.} We note first that
  $$
   0<2-p<\sqrt{2-p}<1
  $$
for $p\in(1,2)$, and
   $$
    1<\sqrt{2-p}<2-p<2
   $$
for $p\in[1/2,1)$. By combining with the monotonicity of $H$ (Lemma \ref{l8.2}), one concludes the non-existence of non-trivial positive classical solution of \eqref{e1.3} for $p\in[1/2,1)\cup(1,2)$ after applying Proposition \ref{p8.1}.

 When $p=1$, solving the linear ordinary differential equation \eqref{e1.3} of second order, one can deduce that all positive classical solutions are given by
    \begin{equation}\label{e8.25}
      h(\theta)=1+A\cos\theta+B\sin\theta, \ \ \forall\theta\in{\mathbb{S}}^1,
    \end{equation}
  where $A, B$ are arbitrary constants satisfying
     \begin{equation}\label{{e8.26}}
       A^2+B^2<1.
     \end{equation}
In the remaining part of (3), the uniqueness result to the logarithmic case $p=0$ has been obtained by Chow \cite{C}. Also, for the centroaffine case $p=-2$, it is well known that by Blaschki-Santalo's inequality, all ellipsoids centered at origin with volume of unit ball are all solutions to \eqref{e1.3} for $p=-2$.

If $p<-7$, applying Lemma \ref{l8.1} and supposing that for some positive integer $\kappa\geq2$ such that
   \begin{equation}\label{{e8.27}}
     2<\kappa\leq\sqrt{2-p},
   \end{equation}
there must be a positive classical solution $u$ of \eqref{e1.3} satisfying
   \begin{equation}\label{e8.28}
      \min_{{\mathbb{S}}^1}u=m, \ \ H(m)=\frac{\pi}{\kappa}.
   \end{equation}
As a result, given each $p<-7$, there exists at least
 $$
   c_p\equiv\big[\sqrt{2-p}\big]_*-1
 $$
positive classical solutions of \eqref{e1.3}, where $[z]_*$ stands for the largest integer no greater than $z$.

Finally, we give the proof of part (4). Noting that for $p<2, p\not=0, 1,-2$, it is inferred from Lemma \ref{l8.3} that there is a small constant $\sigma_p>0$, such that $H$ is a strict monotone function on
   $$
      {\mathcal{\vartheta}}_p\equiv(1-\sigma_p,1+\sigma_p).
   $$
As a corollary,
  \begin{equation}\label{e8.29}
    H(m)\not=\pi/\kappa, \forall m\in{\mathcal{\vartheta}}_p\setminus\{1\}, \ \ \forall \kappa\in{\mathbb{N}}
  \end{equation}
by shrinking the interval ${\mathcal{\vartheta}}_p$ if necessary. Hence part (4) of Theorem \ref{t1.5} follows from Proposition \ref{p8.1}. The proof of Theorem \ref{t1.5} was completed. $\Box$\\

\vspace{40pt}

\section*{Acknowledgments}

The author would like to express his deepest gratitude to Professors Xi-Ping Zhu, Kai-Seng Chou, Xu-Jia Wang and Neil Trudinger for their constant encouragements and warm-hearted helps. This paper is also dedicated to the memory of Professor Dong-Gao Deng.\\


\vspace{40pt}

\begin{thebibliography}{}

\bibitem{A} A.D. Aleksandrov, {\it Existence and uniqueness of a convex surface with a given integral curvature}, C. R. (Doklady) Acad. Sci. URSS (N.S.), {\bf35} (1942), 131-134.\\

\bibitem{An} B. Andrews, {\it Classification of limiting shapes for isotropic curve flows}, J. Amer. Math. Soc., {\bf16} (2003), 443-459.\\

\bibitem{BCD} S. Brendle, K. Choi and P. Daskalopoulos, {\it Asymptotic behavior of flows by powers of the Gaussian curvature}, Acta Math., {\bf219} (2017), 1-16.\\

\bibitem{BLYZ} L.J. B\"{o}r\"{o}czky, E. Lutwak, D. Yang and G. Zhang, {\it The log-Brunn-Minkowski inequality}, Adv. Math., {\bf231} (2012), 1974-1997.\\

\bibitem{C} B. Chow, {\it Deforming convex hypersurfaces by the $n$th root of the Gaussian curvature}, J. Differential Geom., {\bf22} (1985), 117-138.\\

\bibitem{CHLL} S.B. Chen, Y. Huang, Q.R. Li and J.K. Liu, {\it The $L_p$-Brunn-Minkowski inequality for $p<1$}, Advances in Mathematics, {\bf 368} (2020), 107166, 21pp.\\

\bibitem{CL} S.B. Chen and Q.R. Li, {\it On the planar dual Minkowski problem}, Advances in Mathematics, {\bf333} (2018), 87-117.\\

\bibitem{CW} K.S. Chou and X.J. Wang, {\it The $L_p$-Minkowski problem and the Minkowski problem in centroaffine geometry}, Adv. Math., {\bf205} (2006), 33-83.\\

\bibitem{CZ} K.S. Chou and X.P. Zhu, {\it The curve shortening problem}, Chapman $\&$ Hall/CRC, Boca Raton, FL, 2001. x+255 pp. ISBN: 1-58488-213-1.\\

\bibitem{CY} S.Y. Cheng and S.T. Yau, {\it On the regularity of the solution of the $n$-dimensional Minkowski problem}, Comm. Pure Appl. Math., {\bf29} (1976), 495-516.\\

\bibitem{DG} C. Dohmen and Y. Giga, {\it Selfsimilar shrinking curves for anisotropic curvature flow equations}, Proc. Japan Acad. Ser. A, {\bf70} (1994), 252-255.\\

\bibitem{E} L.C. Evans, {\it Partial differential equations}, Second edition. Graduate Studies in Mathematics, {\bf19} American Mathematical Society, Providence, RI, 2010. xxii+749.\\

\bibitem{F} W. Firey, {\it Shapes of worn stones}, Mathematika, {\bf21} (1974), 1-11.\\

\bibitem{G} M.E. Gage, {\it Evolving plane curves by curvature in relative geometries}, Duke Math. J., {\bf72} (1993), 441-466.\\

\bibitem{GLL} P. Guan, J. Li and Y.Y. Li, {\it Hypersurfaces of prescribed curvature measure}, Duke Math. J., {\bf161} (2012), 1927-1942.\\

\bibitem{GM} P. Guan, X.N. Ma, {\it The Christoffel-Minkowski problem. I. Convexity of solutions of a Hessian equation}, Invent. Math., {\bf151} (2003), 553-577.\\

\bibitem{HLX} Y. Huang, J.K. Liu and L. Xu, {\it On the uniqueness of $L_p$-Minkowski problems: The constant $p-$curvature case in ${\mathbb{R}}^3$}, Advances in Math., {\bf 281} (2015), 906-927.\\

\bibitem{HLYZ} Y. Huang, E. Lutwak, D. Yang and G.Y. Zhang, {\it Geometric measures in the dual Brunn-Minkowski theory and their associated Minkowski problems}, Acta Math., {\bf216} (2016), 325-388.\\

\bibitem{J} F. John, {\it Extremum problems with inequalities as subsidiary conditions}, in: Studies and Essays Presented to R. Courant on his 60th Birthday, January 8, 1948, Interscience Publishers, Inc., New York, NY, 1948, pp. 187-204.\\

\bibitem{JLW} H.Y. Jian, J. Lu and X.J. Wang, {\it Nonuniqueness of solutions to the $L_p-$Minkowski problem}, Advences in Mathematics, {\bf281} (2015), 845-856.\\

\bibitem{L} Y.Y. Li, {\it Degree theory for second order nonlinear elliptic operators and its applications}, Comm. Partial Differential Equations, {\bf14} (1989), 1541-1578.\\

\bibitem{L1} E. Lutwak, {\it Dual mixed volumes}, Pacific J. Math., {\bf58} (1975), 531-538.\\

\bibitem{L2} E. Lutwak, {\it The Brunn-Minkowski-Firey theory. I. Mixed volumes and the Minkowski problem}, J. Differential Geom., {\bf38} (1993), 131-150.\\

%\bibitem{LYZ1} E. Lutwak, D. Yang and G. Zhang, {\it $L_p$ affine isoperimetric inequalities}, J. Differential Geom., {\bf56} (2000), 111-132.\\

%\bibitem{LYZ2} E. Lutwak, D. Yang and G. Zhang, {\it Sharp affine $L_p$ Sobolev inequalities}, J. Dffierential Geom., {\bf62} (2002), 17-38.\\

\bibitem{LYZ} E. Lutwak, D. Yang and G. Zhang, {\it On the $L_p-$Minkowski problem}, Trans. Amer. Math. Soc., {\bf356} (2004), 4359-4370.\\

%\bibitem{LYZ4} E. Lutwak, D. Yang and G. Zhang, {\it Optimal Sobolev norms and the $L^p$ Minkowski problem}, Int. Math. Res. Not. IMRN (2006), Art. ID 62987, 21 pp.\\

\bibitem{LO} E. Lutwak and V. Oliker, {\it On the regularity of solution to a generalization of the Minkowski problem}, J. Differential Geom., {\bf40} (1995), 227-246.\\

\bibitem{N} L. Nirenberg, {\it The Weyl and Minkowski problems in differential geometry in the large}, Comm. Pure Appl. Math., {\bf6} (1953), 337-394.\\

\bibitem{N2} L. Nirenberg, {\it Topics in nonlinear functional analyssi}, in: Lecture Notes, 1973-1974, Courant Institute of Mathematical Sciences, New York University, New York, 1974, viii+259 pp.\\

\bibitem{P} A.V. Pogorelov, {\it Extrinsic geometry of convex surfaces},  Translations of Mathematical Monographs, {\bf35}, American Mathematical Society, Providence, RI, 1973, vi+669 pp.

\bibitem{W} H.F. Weinberger, {\it A first course in partial differential equations in complex variables and transform methods}, Blaisdell Publishing Co. Ginn and Co. New York-Toronto-London, 1965 ix+446 pp.\\

\bibitem{Y} H. Yagisita, {\it Non-uniqueness of self-similar shrinking curves for an anisotropic curvature flow}, Cacl. Var. Partial Differential Equations, {\bf26} (2006), 49-55.\\

\end{thebibliography}
\end{document}